\documentclass[12pt]{article}

\font\tenmsb=msbm10
\font\sevenmsb=msbm7
\font\fivemsb=msbm5

\newfam\msbfam
\textfont\msbfam=\tenmsb
\scriptfont\msbfam=\sevenmsb
\scriptscriptfont\msbfam=\fivemsb
\def\Bbb#1{{\fam\msbfam #1}}

\makeatletter
\ifnum\@ptsize=0  \fi \ifnum\@ptsize=2
 \fi 

\newcommand\bR{{\Bbb R}}
\newcommand\bZ{{\Bbb Z}}
\newcommand\bC{{\Bbb C}}
\newcommand\bQ{{\Bbb Q}}

\newcommand\bP{{\Bbb P}}

\newcommand\ra{{\rightarrow}}

\newtheorem{theorem}{Th\' eor\`eme}[section]
\newtheorem{lemma}[theorem]{Lemme}
\newtheorem{corollary}[theorem]{Corollaire}
\newtheorem{proposition}[theorem]{Proposition}
\newtheorem{question}[theorem]{Question}
\newtheorem{re}[theorem]{Remarque}

\newtheorem{definition}[theorem]{Definition}
\newtheorem{conjecture}[theorem]{Conjecture}

\newtheorem{example}[theorem]{Exemple}
\newtheorem{notation}[theorem]{Conventions}

\newenvironment{remark}{\begin{re}\em}{\end{re}}

\begin{document}

\title {FIBRES MULTIPLES SUR LES SURFACES:ASPECTS GEOM\'ETRIQUES, HYPERBOLIQUES ET ARITHM\'ETIQUES} 
\author{Fr\'ed\'eric Campana}

\maketitle

\section*{Introduction}

\

La conjecture de Lang pour les surfaces $S$ de type g\'en\'eral existe en trois versions au moins: hyperbolique, corps de fonctions (complexe) et arithm\'etique. La version hyperbolique faible, par exemple,  affirme l'existence d'une courbe projective $D\subset S$ contenant l'image de toute application holomorphe non-constante de $\Bbb C$ dans $S$. 

Cette conjecture peut-\^etre r\'eduite au cas des courbes, et donc d\'emontr\'ee dans ses trois versions lorsqu'il existe une fibration $f:S\to C$ sur une courbe $C$ de type g\'en\'eral (ie: $g(C)\geq 2$), ou plus g\'en\'eralement, lorsqu'un rev\^etement \'etale fini $S'$ de $S$ poss\`ede cette propri\'et\'e (purement topologique, par un th\'eor\`eme de Siu ([Siu80] ):  le groupe fondamental de $S$ a un sous-groupe d'indice fini admettant pour quotient un groupe de surface). Cette observation est utilis\'ee dans [C-S-S97] et [D97]).

Cette derni\`ere propri\'et\'e peut \^etre cependant aussi formul\'ee comme suit: il existe une fibration $f:S\to C$ dont la {\it base orbifolde} pour les multiplicit\'es {\it classiques} (voir d\'efinition en \ref{bsorb} ci-dessous) est une courbe orbifolde de type g\'en\'eral. La multiplicit\'e {\it classique} d'une fibre $F=\sum_{j\in J} m_j.F_j$  est d\'efinie par $m^*(F):=pgcd\lbrace m_j,j\in J\rbrace$.

On se propose ici d'\'etendre cette solution de la conjecture de Lang au cas o\`u les multiplicit\'es classiques des fibres de $f$ sont remplac\'ees par les multiplicit\'es {\it non classiques}  $m(F):=inf \lbrace m_j,j\in J\rbrace$, introduites naturellement dans [C01] (voir \ref{just} ci-dessous pour une justification). 

Ce changement de d\'efinition soul\`eve des probl\`emes nouveaux: les fibres multiples ne peuvent en effet plus \^etre \'elimin\'ees par un rev\^etement \'etale de $S$, et les solutions de la conjecture de Mordell pour les courbes de type g\'en\'eral (dans ses trois versions) ne peuvent donc plus \^etre invoqu\'ees.

Ce fait est illustr\'e par l'exemple donn\'e au \S \ref{ex} d'une fibration $f:S\to \Bbb P^1$ dans laquelle la surface projective $S$ est lisse et simplement connexe, et la base orbifolde (n\'ecessairement {\it non classique}) est de type g\'en\'eral. De tels exemples peuvent \^etre d\'efinis sur des corps de nombres, ou m\^eme sur $\Bbb Q$.

La solution de la conjecture de Lang dans cette nouvelle situation peut n\'eanmmoins encore \^etre r\'eduite, apr\`es d\'efinition ad\'equate de la notion de point $k$-rationnel d'une courbe orbifolde, \`a la solution d'une version orbifolde de la conjecture de Mordell. Nous \'etablissons ici les deux premi\`eres versions (hyperbolique et corps de fonctions), mais pas la version arithm\'etique. On ne dispose pas, en effet, dans le contexte orbifolde, du plongement Jacobien\footnote{sauf dans le cas d'une courbe elliptique point\'ee. Voir \ref{ell}}, qui joue un r\^ole essentiel dans les solutions existantes de la conjecture de Mordell.

Cette version orbifolde de la conjecture de Mordell arithm\'etique (voir \ref{mo} pour l'\' enonc\'e pr\'ecis), bien que cons\'equence de la conjecture $abc$ est cependant d'un int\'er\^et ind\'ependant, et permettrait par exemple de montrer que les surfaces lisses et simplement connexes construites au \S\ref{ex} ne sont pas potentiellement denses. 

Des discussions avec P. Eyssidieux et E. Peyre sont \`a l'origine de la version corps de fonctions consid\'er\'ee ici. Je les en remercie, ainsi que J.L. Colliot-Th\'el\`ene pour l'observation \ref{abc}, et L.Gruson et T. Peternell, pour leur relecture de la \S.\ref{ex}.

\

\tableofcontents

\

\section{Base orbifolde d'une fibration}

\

On rappelle ici un certain nombre de notions introduites dans [C01] auquel on renvoit pour plus de d\'etails.

\subsection{ Orbifoldes, fibr\'e canonique et dimension de Kodaira}

Une orbifolde $(Y/\Delta)$ est la donn\'ee d'une vari\'et\'e projective complexe $Y$ et d'un $\Bbb Q$-diviseur {\it orbifolde} $\Delta:=\sum_{j\in J} (1-1/m_j).\Delta_j$ sur $Y$, dans lequel $J$ est un ensemble fini, les $m_j>1$ sont des entiers, et les $\Delta_j$ des diviseurs (de Weil) irr\'eductibles distincts de $Y$ tels que le diviseur $K_Y+\Delta$ soit $\Bbb Q$-Cartier. (Cette condition sera toujours v\'erifi\'ee dans la suite o\`u $Y$ sera lisse). On dira que l'orbifolde $(Y/\Delta)$ est {\it support\'ee} par $Y$.

Le $\Bbb Q$-diviseur $K_Y+\Delta:=K_{(Y/\Delta)}$ est appel\'e le {\it diviseur canonique} de $(Y/\Delta)$, et $\kappa(Y/\Delta):=\kappa(Y,K_Y+\Delta)$ est appel\'e sa dimension de Kodaira. On a donc toujours: $\kappa(Y/\Delta)\geq \kappa(Y)$ si $Y$ est lisse (ou si $K_Y$ est $\Bbb Q$-Cartier).

On dira que $(Y/\Delta)$ est {\it de type g\'en\'eral} si $\kappa(Y/\Delta)=dim(Y)>0$.

On pourrait, dans cette situation, d\'efinir aussi le groupe fondamental de cette orbifolde. Voir [C01].

\subsection{Base orbifolde d'une fibration} \label{bsorb}

\subsubsection{Multiplicit\'es}

Soit $F:=\sum_{k\in K}m_k.F_k$ un cycle analytique de dimension pure $p\geq 0$ d'un espace analytique $Z$. Ceci signifie que les $F_k$ sont des sous-ensembles analytiques compacts irr\'eductibles distincts de $Z$ de dimension complexe $p$, $K$ un ensemble fini, et les $m_k>0$ des entiers. La {\it multiplicit\'e} $m(F)$ (resp. la {\it multiplicit\'e classique}) $m^*(F)$) est le plus petit (resp. le plus grand commun diviseur) des $m_k,k\in K$. Donc $m^*(F)$ divise $m(F)$.

\

\begin{remark} \label{just} L'introduction de la notion non classique de multiplicit\'e provient de ce qu'on a une bijection naturelle entre les fibrations de type g\'en\'eral au sens non classique et les faisceaux de Bogomolov; de plus, elle est mieux adapt\'ee \`a l'etude de la pseudom\'etrique de Kobayashi (et probablement aussi \`a la g\'eom\'etrie arithm\'etique). Voir [C01] pour une justification plus d\'etaill\'ee.
\end{remark}

\subsubsection{Base orbifolde d'une fibration}

Une {\it fibration} d\'esigne dans le pr\'esent texte une application holomorphe surjective et \`a fibres connexes $f:X\ra Y$ entre vari\'et\'es projectives complexes lisses et connexes.

Si $f:X\to Y$ est une fibration, et si $D$ est un diviseur irr\'eductible de $Y$, on notera $f^*(D)=F+R$, o\`u $F$ (resp. $R$) est la r\'eunion des composantes irr\'eductibles de $f^*(D)$ dont l'image par $f$ est $D$ (resp. est de codimension au moins $2$ dans $Y$). On \' ecrit $F:=\sum_{k\in K}m_k.F_k$, et on d\'efinit alors la multiplicit\'e (resp. la {\it multiplicit\'e classique} de la fibre g\'en\'erique de $f$ au-dessus de $D$ comme \'etant: $m(f,D):=m(F)$ (resp. $m^*(f,D):=m^*(F)$).

Donc $m(f,D)$ et $m^*(f,D)$ sont \'egaux \`a $1$ sauf pour un nombre fini de $D$, contenus dans le lieu au-dessus duquel $f$ n'est pas lisse.

On d\'efinit alors la {\it base orbifolde} (resp. la {\it base orbifolde classique}) de $f$ comme \'etant le couple $(Y/\Delta(f))$ (resp. $(Y/\Delta^*(f))$), avec: 

$\Delta(f):=\sum_{D\subset Y}(1-(1/m(f,D))).D$ (resp.  $\Delta^*(f):=\sum_{D\subset Y}(1-(1/m^*(f,D))).D$.

On notera que ces sommes sont bien finies, puisque $m(f,D)$ (et a fortiori $m^*(f,D)$) sont \'egales \`a $1$ sauf pour un nombre fini de $D$.

On a alors aussi: $\kappa(Y)\leq \kappa(Y/\Delta^*(f))\leq \kappa(Y/\Delta(f))$.

On dira que $f$ est {\it de type g\'en\'eral} (resp. {\it de type g\'en\'eral au sens classique}) si $Y$ est une courbe\footnote{Si $Y$ est de dimension $2$ ou plus (cas non utilis\'e ici), la d\'efinition est plus compliqu\'ee. Voir [C01].} et si la base orbifolde $(Y/\Delta(f))$ (resp. la base orbifolde classique $(Y/\Delta^*(f))$ de $f$ est de type g\'en\'eral (ie: de dimension de Kodaira $1$).

\

\begin{remark} Dans les situations consid\'er\'ees ici, $X$ sera une surface, et $Y$ une courbe. Dans ce cas, les diviseurs $D$ de $Y$ sont des points. De plus, il r\'esulte alors de [B-P-V,V. 7. p.150]  que si $f$ est une fibration elliptique, alors $m^*(f,D)=m(f,D)$ pour tout $D\subset Y$. Si la fibre g\'en\'erique de $f$ est rationnelle, alors $m(f,D)=m^*(f,D)=1$ pour tout $D\subset Y$. La distinction entre $\Delta(f)$ et $\Delta^*(f)$ n'apparait donc que si les fibres lisses de $f$ sont des courbes de type g\'en\'eral, condition  satisfaite si $X$ est elle-m\^eme de type g\'en\'eral.
\end{remark}

\subsection{Courbes orbifoldes, morphismes d'orbifoldes}\label{corb}

\

On va consid\'erer ici une orbifolde $(C/\Delta)$, dans laquelle $C$ est une courbe projective lisse et connexe de genre $g(C)$ , et $\Delta:=\sum_{k=1}^{k=N}(1-1/m_k). p_k$, o\`u les $p_k$ sont des points (donc des diviseurs irr\'eductibles) de $C$. 

On notera parfois, lorsque la connaissance des points $p_k$ est superflue, $(C/(m_1,...,m_N))$ une telle orbifolde, les $m_k$ \'etant ordonn\'es de telle sorte que $2\leq m_1\leq m_2\leq ...\leq m_N$.

Par exemple, $\Bbb P^1/(2,3,7)$ d\'esigne une orbifolde $\Bbb P^1/\Delta$, avec $\Delta= (1/2).p_1+(2/3).p_2+(6/7). p_3$, o\`u les $p_j\in \Bbb P^1, j=1,2,3$ sont distincts.

Si $(C/\Delta)$ et $(C'/\Delta')$ sont des orbifoldes dont les supports sont des courbes, un {\it morphisme  d'orbifoldes} $g:(C/\Delta)\to (C'\Delta')$ est la donn\'ee d'une application holomorphe surjective $g:C\to C'$ telle que $g^*(K_{C'/\Delta'})\subset K_{C/\Delta}$. On dit que $g$ est {\it \'etale} si l'on a l'\'egalit\'e $g^*(K_{C'/\Delta'})=K_{C/\Delta}$.
Alors $\kappa(C/\Delta)\geq \kappa(C'\Delta')$, avec \'egalit\'e si $g$ est \'etale.

Donc $g:(C/\Delta)\to (C'/\Delta')$ est un morphisme d'orbifoldes si et seulement si, pour tout $x\in C$, on a: $r(g,x).\Delta(x)\geq \Delta(g(x))$, o\`u $r(g,x)$ est l'indice de ramification de $g$ en $x$, $\Delta(x)$ est la multiplicit\'e de $\Delta$ en $x$, \'egale \`a $1$ si $x$ n'est pas l'un des $p_k$, et \'egale \`a $m_k$ si $x=p_k$. On d\'efinit $\Delta'(x')$ de mani\`ere similaire si $x'\in C'$. Et $g$ est \' etale exactement lorsque toutes ces in\'egalit\'es sont des \'egalit\'es.

\begin{example}\label{etale}

\

1. Soit $g:C:=E\to \Bbb P^1:=C'$ un rev\^etement double ramifi\'e en $4$ points $a_k\in \Bbb P^1, k=1,2,3,4$, avec $E$ une courbe elliptique. Soit $\Delta$ le diviseur vide sur $E$ et $\Delta':=\sum_{k=1}^{k=4}(1/2).a_k$. Alors $g:(E/\emptyset)\to (\Bbb P^1/\Delta')$ est un morphisme d'orbifoldes \'etale de degr\'e $2$, ainsi que $g:(E/((1/2).\bar a_1))\to (\Bbb P^1/\Delta")$, si $\bar a_1$ est l'unique point de $E$ dont l'image par $g$ est $a_1$, et si $\Delta":=(3/4). a_1+\sum_{k=2}^{k=4}(1/2).a_k$.

2. Plus g\'en\'eralement, si $(C':=\Bbb P^1/\Delta')$ est une orbifolde telle que $N\geq 3$ dans la notation pr\'ec\'edente, il existe un morphisme d'orbifolde \'etale $g:C\to (\Bbb P^1/\Delta')$. (On note $C$ une orbifolde support\'ee par une courbe et de diviseur vide). Voir, par exemple, [N87, p. 26] pour cette assertion. Donc: $\kappa(\Bbb P^1/\Delta)=1$ (resp. $0$) si et seulement s'il existe un morphisme \' etale $g:C\to (\Bbb P^1/\Delta)$, avec $g(C)\geq 2$ (resp. $g(C)=1$).

\end{example}

\begin{remark}\label{comp}
S'il existe un morphisme d'orbifolde $g:(C/\Delta)\to (C'\Delta')$, on note: $(C/\Delta)\geq (C'\Delta')$.
On v\'erifie que $\kappa(\Bbb P^1/\Delta)=0$ exactement dans les cas suivants: $N=3$ et $\Delta=(3,3,3)$, $\Delta=(2,3,6)$ ou $\Delta=(2,4,4)$; ou bien $N=4$ et $\Delta=(2,2,2,2)$.
De plus: $\kappa(\Bbb P^1/\Delta)=-\infty$ si et seulement si l'on est dans l'un des cas suivants: $N\leq 2$, $N=3$ et $\Delta$ est l'une des suivantes: $(2,2,m)$, $(2,3,m')$, avec $m$ arbitraire, $m'\leq 5$.

On en d\'eduit que $\kappa(\Bbb P^1/\Delta)=1$ si et seulement si $(\Bbb P^1/\Delta)\geq (\Bbb P^1/\Delta')$, o\`u $\Delta'$ est l'une des cinq suivantes: $(2,3,7), (2,4,5), (3,3,4), (2,2,2,3), (2,2,2,2,2)$.

\end{remark}

\

\subsection{Fibrations de type g\'en\'eral sur une surface projective}

\

On consid\`ere dans cette section une fibration $f:X\to C$ dans laquelle $X$ est une {\it surface} (projective complexe). On note $F$ une fibre lisse de $f$, $g(F)$ son genre, et $\kappa(X)\in \lbrace -\infty, 0,1,2\rbrace$ la dimension canonique (ou de Kodaira) de $X$.

\begin{proposition} Si $f$ est de type g\'en\'eral, alors:

1. $\kappa(X)\neq 0$. De plus:

2. Si $\kappa(X)=-\infty$, alors: $g(F)=0$ et $g(C)=q(X)\geq 2$. De plus, $f$ est \`a la fois le quotient rationnel, le morphisme d'Albanese, et le coeur de $X$ (voir [C01] pour ce terme).

3. Si $\kappa(X)=1$, alors $g(F)=1$, et $f$ est \`a la fois la fibration d'Iitaka-Moishezon  et le coeur de $X$. De plus, $f$ est de type g\'en\'eral au sens classique, donc il existe un rev\^etement \'etale fini $u:X'\to X$ tel que $g(C')\geq 2$ si $f':X'\to C'$ est la  {\it partie connexe} de la factorisation de Stein de $v\circ f'=f\circ u:X'\to C$, et $v:C'\to C$ sa {\it partie finie}. 

4. $\kappa(X)=2$ si et seulement si $g(F)\geq 2$.
 
\end{proposition}

\

{\bf D\'emonstration de 2.1:} Puisque $f$ est une fibration de type g\'en\'eral, le th\'eor\`eme d'additivit\'e orbifolde (voir [C01], \S 4) affirme que $\kappa(X)=\kappa(F)+dim_{\bC}(C)=\kappa(F)+1$. Ceci montre donc que $\kappa(F)=-\infty$ si $\kappa(X)=-\infty$, que $\kappa(X)\neq 0$, et que $\kappa(F)=\kappa(X)-1$ si $\kappa(X)=1,2$. Si $\kappa(X)\neq 0,2$, on voit donc que $f$ est \`a la fois \`a fibres {\it sp\'eciales} (voir [C01]) et de type g\'en\'eral (comme fibration). C'est donc le {\it coeur} de $X$ (voir [C01] pour la description explicite du coeur dans le cas des surfaces). Cette description donne les autres assertions de l'\'enonc\'e.

\

\begin{remark} Le r\'esultat pr\'ec\'edent est un cas tr\`es particulier de r\'esultats valables en dimension sup\'erieure. Voir [C01].
\end{remark}

\begin{remark} Le fait d'\^etre de type g\'en\'eral au sens classique s'interpr\`ete g\'eom\'etriquement tr\`es simplement (voir [Ca01]): $f$ est de type g\'en\'eral au sens classique si et seulement si, apr\`es changement de base fini ad\'equat (ramifi\'e en g\'en\'eral) $v:C'\to C$ et normalisation, le morphisme $f':X'\to C'$ d\'eduit de $f$ n'a pas de fibre multiple au sens classique, et si $C'$ est une courbe de type g\'en\'eral. Donc si $f:X\to C$ est une fibration de type g\'en\'eral, le groupe fondamental de $X$ admet un sous-groupe d'indice fini dont un quotient est le groupe fondamental d'une courbe de genre $2$ ou plus. (La r\'eciproque est d'ailleurs vraie aussi, par un th\'eor\`eme de Siu ([Siu80])). 

Nous verrons par contre l'existence de fibrations de type g\'en\'eral $f:X\to C$ avec $X$ simplement connexe, et donc $C\cong \Bbb P^1$. Une telle $f$ n'est donc pas de type g\'en\'eral au sens classique.
\end{remark}

\

\begin{remark} \label{reduc} Si $f:X\to C$ est une fibration de type g\'en\'eral {\it au sens classique} avec $X$ une surface, les conjectures arithm\'etiques et hyperboliques de Lang peuvent \^etre r\'eduites aux \'enonc\'es analogues pour les courbes, et donc r\'esolues, \`a l'aide de r\'esultats bien connus. 

Expliquons le principe de cette r\'eduction dans le cas arithm\'etique (les autres cas sont similaires). Elle se fait en les quatre \'etapes ci-dessous, supposant $u:X'\to X$ \'etale et d\'efinie sur le corps de nombres $k$, ainsi que $f':X'\to Y'$, o\`u $Y'$ est une courbe de type g\'en\'eral.

0. Par le th\'eor\`eme de Chevalley-Weil, il existe une extension finie de corps $k'/k$ telle que $u(X'(k'))\supset X(k)$. Il suffit donc d'\'etablir le r\'esultat pour $X'$.

1. $f'(X'(k'))\subset Y'(k')$.

2. $Y'(k')$ est fini (conjecture de Mordell, \'etablie dans [F83]). C'est l'\'etape cruciale.

3. Si $V'\subset Y'$ est l'ouvert au-dessus duquel $f'$ est lisse, et si $U':=(f')^{-1}(V')$, alors les fibres de $f':X'(k')\cap U'\to Y'(k')\cap V'$ sont finies, par Mordell encore.

Par contre, lorsque $f:X\to C$ est de type g\'en\'eral au sens non classique ci-dessus, des arguments nouveaux sont n\'ecessaires, d\^us au fait que l'\'etape 0. ci-dessus ne peut plus \^etre appliqu\'ee. On va \'etablir dans cette situation les versions hyperboliques et corps de fonctions. La version arithm\'etique conduit \`a une version orbifolde de la conjecture de Mordell (voir \ref{mo'}) pour laquelle les d\'emonstrations connues ne semblent pas s'adapter, car bas\'ees de mani\`ere essentielle sur l'\'etude de la Jacobienne.
\end{remark}

\section{Version hyperbolique}

\subsection{Pseudom\'etrique de Kobayashi}

On consid\`ere dans cette section aussi une fibration $f:X\to C$ dans laquelle $X$ est une {\it surface projective complexe}. On note encore $F$ une fibre lisse de $f$, et $\kappa(X)\in \lbrace -\infty, 0,1,2\rbrace$ la dimension canonique (ou: de Kodaira) de $X$. On note $d_X:X\times X\to \bR^+$ la pseudom\'etrique de Kobayashi (voir [K71]). Si $f:X\to C$ est une application, et $\delta_C$ une pseudom\'etrique sur $C$, on note $f^*(\delta_C)$ la pseudom\'etrique sur $X$ d\'efinie par: $f^*(\delta_C):=(\delta_C)\circ(f\times f):X\times X\to \bR^+$. Si $d,d'$ sont deux pseudom\'etriques sur $X$, on d\'efinit l'in\'egalit\'e $d'\geq d$ de la mani\`ere \'evidente. Une pseudom\'etrique $d$ est une m\'etrique si $d(x,y)>0$ lorsque $x\neq y$. En g\'en\'eral, $d_X$ n'est pas une m\'etrique; par exemple $d_{\bC}\equiv 0$.

Rappelons que $d_X$ est la plus grande des pseudom\'etriques sur $X$ telle que $d_X\leq h^*(d_{D})$, si $d_{D}$ est la m\'etrique de Poincar\'e sur le disque unit\'e $D$ de $\bC$, et $h:D\to X$ une application holomorphe arbitraire. (Le lemme d'Ahlfors-Schwartz montre que les deux sens de la notation $d_D$ coincident). Il est \'evident que si $g:Y\to X$ est une application holomorphe arbitraire, on a: $d_Y\geq g^*(d_X)$.  (Propri\'et\'e de d\'ecroissance de la pseudom\'etrique de Kobayashi  par applications holomorphes).

\begin{remark} Si $f:X\to C$ est de type g\'en\'eral au sens classique, apr\`es rev\^etement \'etale fini ad\'equat  $u:X'\to X$, la {\it partie connexe} $f':X'\to C'$ de la factorisation de Stein de $f\circ u$ a pour image $C'$, une courbe de type g\'en\'eral. On en d\'eduit ais\'ement que:
 $d_X\geq f^*(\delta_C)$, o\`u $\delta_C$ est une m\'etrique sur $C$ qui d\'efinit la topologie {\it analytique} de $C$. En particulier, si $X$ (et donc la fibre g\'en\'erique $F$ de $f$) est de type g\'en\'eral, alors $d_X$ est une m\'etrique sur l'ouvert de Zariski $U$ de $X$ qui est la r\'eunion des fibres lisses de $f$. (Les arguments sont des cons\'equences imm\'ediates du th\'eor\`eme de Liouville et du fait que les courbes de genre $2$ ou plus sont uniformis\'ees par le disque unit\'e).
 \end{remark}

\subsection{Pseudom\'etrique de Kobayashi orbifolde}

Nous allons, suivant [C01], \'etendre la notion de pseudom\'etrique de Kobayashi au cadre orbifolde:

\begin{definition} Soit $(C/\Delta)$ une orbifolde, o\`u $C$ est une courbe projective complexe, et $\Delta:=\sum_{j\in J} (1-1/m_j).p_j$ un $\bQ$-diviseur, les $p_j$ \'etant des points distincts de $C$ et $J$ un ensemble fini. On d\'efinit $d_{(C/\Delta)}$ comme \'etant la plus grande des pseudom\'etriques $\delta$ sur $C$ telles que $\delta\leq h^*(d_D)$, ceci pour toutes les $h:D\to C$, applications holomorphes {\it compatibles avec le diviseur orbifolde $\Delta$}, c'est-\`a-dire telles que $h^*(p_j)\geq m_j.h^{-1}(p_j)$, pour tout $j\in J$. (Autrement dit: $h(z)$ ramifie au moins \`a l'ordre $m_j$ en tout $z\in D$ tel que $h(z)=p_j$).
\end{definition}

\

\begin{remark} Si $\Delta=\emptyset$, on a donc: $d_{(C/\emptyset)}=d_C$. De plus, si $g:C'\to (C/\Delta)$ est une application holomorphe de la courbe $C'$ vers $C$ compatible avec $\Delta$, alors $d_{C'}\geq g^*(d_{(C/\Delta)})$, puisque toute application $h:D\to C'$ fournit par composition avec $f$ une application $f\circ h:D\to C$ compatible avec le diviseur orbifolde $\Delta$. Le m\^eme argument fournit, plus g\'en\'eralement le r\'esultat suivant (que l'on pourrait d'ailleurs aussi formuler lorsque $X$ aussi est munie d'un diviseur orbifolde):
\end{remark}

\
\begin{proposition}\label{kobcorb}

\

 1. Soit $f:X\to C$ une fibration sur une courbe. Alors: $d_X\geq f^*(d_{(C/\Delta(f))})$.

2. Soit $g:(C/\Delta)\to (C'/\Delta')$ un morphisme d'orbifoldes support\'ees par des courbes $C$ et $C'$. Alors $d_{(C/\Delta)}\geq g^*(d_{(C'/\Delta')})$. En particulier, $d_{(C/\Delta)}$ est une m\'etrique sur $C$ si $d_{(C'/\Delta')}$ est une m\'etrique sur $C'$.

\end{proposition}

\

On \'etablit dans [C-W04] le r\'esultat suivant, comme cons\'equence d'une version orbifolde du th\'eor\`eme de Brody:

\

\begin{proposition} Soit $(C/\Delta)$ une courbe projective complexe connexe munie d'une structure d'orbifolde $\Delta$. Alors $d_{(C/\Delta)}$ est une m\'etrique sur $C$ si et seulement si $(C/\Delta)$ est de type g\'en\'eral. Sinon, $d_{(C/\Delta)}$ est nulle.
\end{proposition}

\begin{remark} La d\'efinition donn\'ee de la pseudom\'etrique de Kobayashi orbifolde, o\`u l'ordre de ramification de $h$ en $z$ tel que $h(z)=p_j$ est seulement suppos\'e \^etre au moins \'egal \`a $m_j$ (mais n'est pas suppos\'e \^etre un multiple de $m_j$) interdit de d\'emontrer le r\'esultat pr\'ec\'edent en choisissant un rev\^etement ramifi\'e $g:C'\to C$ de $C$ ramifiant \`a l'ordre $m_j$ en chaque point $p'\in C'$ tel que $g(p')=p_j$, et en relevant \`a $C'$ les applications holomorphes $h:D\to C$ compatibles avec $\Delta$.

Utilisant les notations introduites dans \ref{corb} on voit donc que les orbifoldes de type $(2,2,2,2)$, $(3,3,3)$, $(2,3,6)$ sur $\bP^1$ ont une pseudom\'etrique de Kobayashi nulle. (Remarquons q'elles ont d'ailleurs des rev\^etements cycliques de degr\'es respectifs $2,3,6$, \'etales au sens orbifolde qui sont des courbes elliptiques, d\'efinis (dans l'espace total des fibr\'es $\cal O$$(4)$,$\cal O$$(3)$ et $\cal O$$(6)$) par les \'equations $z^2=x(x-1)(x-\lambda)$, $z^3=x(x-1)$, et $z^6=x^3.(x-1)^2$, $z$ (resp. $x$) \'etant la coordonn\'ee de la fibre (resp. de la base).
\end{remark}

\

\subsection{Non d\'eg\'en\'erescence de la pseudom\'etrique de Kobayashi}

\

\begin{proposition} {\it Soit $f:X\to C$ une fibration de type g\'en\'eral de la surface $X$ sur la courbe $C$.

1. Si $\kappa(X)=-\infty$,  $d_X=f^*(d_C)$.

2. Si $\kappa(X)=1$, alors $ d_X=f^*(\delta _C)$, o\`u $\delta _C$ est une pseudom\'etrique sur $C$ telle que $\delta _C \geq f^*(d_{(C/\Delta(f))})$, $d_{(C/\Delta(f))}$ \'etant la pseudom\'etrique de l'orbifolde $(C/\Delta(f))$, d\'efinie ci-dessus, et qui est une m\'etrique dans ce cas.

3. Si $\kappa(X)=2$, alors $d_X$ est une m\'etrique sur l'ouvert de Zariski $U$ de $X$ constitu\'e de la r\'eunion des fibres lisses de $f$.}
\end{proposition}

\

{\bf D\'emonstration:} La d\'emonstration de la premi\`ere assertion peut \^etre faite sans recourir \`a 3.4 ci-dessus: si $\kappa(X)=-\infty$, alors $X$ est birationnelle \`a $\bP_1\times C$, d'o\`u la conclusion, puisque la pseudom\'etrique est un invariant birationnel. (On peut aussi se dispenser de 3.4 pour la seconde assertion, en consid\'erant des changements de base fini ad\'equats).

D\'emontrons l'assertion 2. Puisque $d_X$ est nulle sur la fibre g\'en\'erale (elliptique) de $f$, elle est nulle sur toute fibre de $f$, puisque continue pour la topologie {\it analytique}. Il existe donc une unique pseudom\'etrique $\delta _C$ telle que $d_X=f^*(\delta _C)$. L'in\'egalit\'e $d_X\geq f^*(d_{(C/\Delta(f))})$ r\'esulte de 3.3, et montre que $\delta _C\geq d_{(C/\Delta(f))}$. D'apr\`es 3.4, $d_{(C/\Delta(f))})$ est une m\'etrique sur $C$. D'o\`u les assertions. 

D\'emontrons l'assertion 3. On a, par 3.3 et 3.4: $d_X\geq f^*(d_{(C/\Delta(f))})$, et $\delta_C:=d_{(C/\Delta(f))})$ est une m\'etrique sur $C$. Soit $a$  un point de $U$. Si $b\in X$ est tel que $d_X(a,b)<2.\epsilon$, o\`u  $\epsilon >0$ est donn\'e, il existe une {\it D-chaine} de $X$ de longueur inf\'erieure \`a $\epsilon$ joignant $a$ et $b$. C'est-\`a-dire qu' il existe une suite d'applications holomorphes $h_i:D\to X$ et de points $a_i,i=0,1,...,n+1$ tels que $a_0=a, a_{n+1}=b$ et de $z_i\in D$ tels que $h_i(0)=a_i$, $h_i(z_i)=a_{i+1}$ et $\sum_{i=0}^{i=n}\vert z_i\vert <\epsilon$. Puisque $d_X\geq f^*(\delta_C)$, tous les $a_i$ sont dans la boule ouverte $B_{\epsilon}$ dans $C$, de centre $c:=f(a)$ et de rayon $\epsilon$ pour la distance $\delta_C$. Comme $\delta_C$ d\'efinit sur $C$ la topologie analytique, $B_{\epsilon}$ est contenue dans un disque analytique $B$ centr\'e en $c$, et tel que $V:=f^{_1}(B)\subset U$, si $\epsilon$ est assez petit. Mais alors $d_{X\vert V}=d_V$ est une m\'etrique, par [K98,3.11.2], puisque $F:=f^{-1}(c)$ est hyperbolique.

Nous allons maintenant d\'emontrer autrement l'assertion 3 pr\'ec\'edente sans recourir \`a 3.4. En fait, c'est une cons\'equence imm\'ediate du lemme suivant:

\ 

\begin{lemma} Soit $f:X\to C$ une fibration de type g\'en\'eral sur une courbe $C$ (la vari\'et\'e $X$ est de dimension arbitraire). 

1. $f\circ h:\bC\to C$ est constante, pour toute application holomorphe $h:\bC\to X$.

2. Soit $\delta _C$ la pseudom\'etrique sur $C$ telle que $\delta_C(a',b')=inf \lbrace d_X(a,b)\vert f(a)=a',f(b)=b'\rbrace, \forall a',b' \in C$. Alors $\delta _C$ est une m\'etrique sur $C$ si $X$ est de type g\'en\'eral.
\end{lemma}

\

{\bf D\'emonstration:} La premi\`ere assertion est le cas particulier $p=1$ du th\'eor\`eme 8.1 de [C01], qui affirme (entre autres) que si $g:X\to Y$ est une fibration de type g\'en\'eral avec $dim(Y)=p>0$, alors $g\circ h:\bC^p\to Y$ est d\'eg\'en\'er\'ee, pour toute application m\'eromorphe $h:\bC^p\to X$.

Pour \'etablir l'assertion  $2$, on utilise [Ko98, thm 3.5.31,p. 94] qui affirme que $d_X(a,b)=inf_{c} \lbrace \int_cF_X(c')\rbrace$, si $c:[0,1]\to X$ est un chemin de classe $C^1$ par morceaux joignant $a$ et $b$ et de d\'eriv\'ee $c'$, et o\`u $F_X:TX\to [0,+\infty[$ est la pseudom\'etrique de Kobayashi infinit\'esimale sur le fibr\'e tangent $TX$ de $X$.

Supposons qu'il existe $a,b\in X$ tels que $d_X(a,b)=0$ et $f(a)\neq f(b)$. On peut supposer que $a$ et $b$ sont dans l'ouvert $U$.
On va montrer qu'il existe un vecteur tangent $t\in TX$ tel que $f_*(t)\neq 0$, et tel que $F_X(t)=0$. 

En effet, sinon $F_X(t)\geq A. \vert f_*(t)\vert$, pour une constante $A>0$, si $\vert t'\vert$ est une m\'etrique hermitienne continue sur $TC$, et si le point base $m$ de $t\in T_mX$ est dans un ouvert connexe $U':=f^{-1}(V')$ de $X$ contenant $a$ et $b$ et relativement compact dans $U$ . Ceci r\'esulte de la continuit\'e de $F_X$, \'etablie dans [Wr77] pour les vari\'et\'es de type g\'en\'eral. 

On en d\'eduit donc que $d_X(a,b)\geq A. \int _{c\cap V} \vert f_*(c')\vert \geq A. dist_V(a',b')>0$, avec: $a':=f(a),b':=f(b)$, pour tout chemin $c$, si $dist$ est la distance sur $C$ d\'eduite de $\vert t'\vert $ sur $TC$, et si $dist_V(a,b)$ est la somme des distances (au sens de $dist$) de $a'$ \`a $C/V'$ et de $b'$ \`a $C/V'$. On en d\'eduit donc une contradiction: il existe bien $t\in TU$ avec les propri\'et\'es annonc\'ees.

Le lemme de reparam\'etrisation de Brody montre alors qu'il existe une application holomorphe $h:\bC\to X$ telle que $h'(0)=t$. Ceci contredit l'assertion 1. Donc un tel couple $(a,b)$ n'existe pas, et $\delta _C$ est bien une m\'etrique sur $C$.

\
\begin{question} A-t'on: $\delta_C=d_{(C/\Delta(f))})$ dans l'assertion 2 de la proposition pr\'ec\'edente?
Peut-on donner une description explicite des m\'etriques de Kobayashi sur les orbifoldes de type $(m_1,m_2,m_3)$ sur $\bP^1$? Par exemple, ces m\'etriques coincident-elles avec celles de leurs rev\^etements orbifoldes \'etales, qui sont des courbes de genre $2$ ou plus sans structure orbifolde?
\end{question}

\

\begin{remark} On peut, plus g\'en\'eralement, suivant [C01], d\'efinir la pseudom\'etrique de Kobayashi d'une orbifolde $(X/\Delta)$ de dimension arbitraire (o\`u $\Delta:=\sum_{j\in J} (1-1/m_j).D_j$ est un $\bQ$-diviseur effectif sur $X$ , les $D_j$ \'etant des diviseurs irr\'eductibles r\'eduits de $X$ et $J$ un ensemble fini. On d\'efinit $d_{(X/\Delta)}$ comme \'etant la plus grande des pseudom\'etriques $\delta$ sur $X$ telles que $\delta\leq h^*(d_D)$, ceci pour toutes les $h:D\to X$, applications holomorphes {\it compatibles avec le diviseur orbifolde $\Delta$}, c'est-\`a-dire telles que $h^*(D_j)\geq m_j.h^{-1}(D_j)$, pour tout $j\in J$. (Autrement dit: $h(z)$ ramifie au moins \`a l'ordre $m_j$ en tout $z\in D$ tel que $h(z)\in D^0_j$, si $D^0_j$ est l'ouvert de $D_j$ constitu\'e des points n'appartenant pas \`a l'une des autres composantes de $\Delta$).

Les propri\'et\'es \'enonc\'ees ci-dessus lorsque $X$ est une courbe restent valables en dimension sup\'erieure, par les m\^emes arguments.
\end{remark}

\section{Version Corps de Fonctions}

On traduit en termes g\'eom\'etriques les notions de g\'eom\'etrie alg\'ebrique usuelles lorsque le corps de base est le corps $k_B$ des fonctions m\'eromorphes sur une courbe projective complexe lisse et connexe $B$. On trouvera dans [C01] d'autres remarques sur cette situation. Voir [Ca 04] pour les probl\`emes d'effectivit\'e et d'uniformit\'e.

\subsection{Vari\'et\'es d\'efinies sur un corps de fontions complexes}\label{varcf}

Soit $B$ une courbe projective complexe lisse et connexe, $k_B$ d\'esignant son corps des fonctions m\'eromorphes. On a donc correspondance bijective entre les extensions finies de $k_B$ et les rev\^etements ramifi\'es $w:B'\to B$.

Un vari\'et\'e $X_B$ de dimension $n$ sur $k_B$ est une vari\'et\'e projective complexe $X$ de dimension $n+1$ munie d'une fibration $g:X\to B$. (Le {\it mod\`ele } $(X,g)$ est en fait d\'efini \`a \'equivalence birationnelle pr\`es, nous ne consid\'erons ici que des propri\'et\'es birationnelles, et nous ferons donc des changements de mod\`ele sans le pr\'eciser). Une fibration $f_B:X_B\to Y_B$ sera simplement une fibration $f:X\to Y$ telle que $h\circ f=g$, si $Y_B=(Y,h)$, avec $h:Y\to B$ la fibration d\'efinissant la $k_B$-vari\'et\'e $Y_B$.

Un point $k_B$-rationnel de $X_B$ est une section $s:B\to X$ \`a $g$. On notera $X(k_B)$ l'ensemble de ces sections. La conjecture de Lang affirme que cet ensemble est fini si $X_B$ est de type g\'en\'eral (ie: si sa fibre g\'en\'erique $X_b:=g^{-1}(b)$ l'est), et si $X_B$ n'est pas isotriviale). Voir [C01] pour des r\'ef\'erences classiques et des observations sur cette question.

Ces notions sont compatibles avec les changements de base finis (ou extensions de corps de base finies) $w:B'\to B$. En particulier, on a une inclusion naturelle $X(k_B)\subset X(k_{B'})$ si $k_B\subset k_{B'}$, et une application naturelle $f_*:X(k_B)\to Y(k_B)$ si $f_B:X_B\to Y_B$ est une fibration de $k_B$-vari\'et\'es.

Remarquons que si $Y_B$ est une courbe sur $k_B$, et si $Y'_B$ est une courbe birationnelle \`a $Y_B$ (ie: on a une application birationnelle $h$ au-dessus de $B$ entre les surfaces $Y$ et $Y'$), alors $Y(k_B)=Y'(k_B)$. En dimension sup\'erieure, on a \'egalit\'e pour les points $k_B$-rationnels ext\'erieurs aux diviseurs exceptionnels de $h$.

\subsection{Orbifoldes d\'efinies sur un corps de fonctions}

Une orbifolde $(X_B/\Delta_B)$ sur $k_B$ est donc la donn\'ee d'une vari\'et\'e projective $g:X\to B$ sur $k_B$ et d'un diviseur orbifolde $\Delta_B=\sum_{j\in J}(1-1/m_j).\Delta_j$ de $X_B$. Les $\Delta_j$ sont des diviseurs irr\'eductibles distincts de $X$ . Par intersection avec la fibre g\'en\'erique $X_b:=g^*(b)$ de $g$, $\Delta_B$ induit une orbifolde $(X_b/\Delta_b)$ de $X_b$, ayant pos\'e:  $\Delta_b=\sum_{j\in J}(1-1/m_j).\Delta_{j,b}$, et $\Delta_{j,b}:=\Delta_j\cap X_b$.

La dimension de Kodaira $\kappa(X_B/\Delta_B)$ n'est autre que $\kappa(X_b/\Delta_b)$, pour $b\in B$ g\'en\'eral.

\begin{example} L'exemple crucial d'orbifolde sur $k_B$ consid\'er\'e ici est bien s\^ur la base orbifolde $Y/\Delta(f)$ d'une fibration $f:X\to Y$ entre $B$-vari\'et\'es $X_B$ et $Y_B$ d\'efinies par des fibrations $g:X\to B$ et $h:Y\to B$. La condition pour que $f$ soit d¥\'efinie sur $k_B$ est que $f$ soit au-dessus de $B$, c'est-\`a-dire que $g=f\circ h$.

L'orbifolde $(Y_b/\Delta(f)_b)$ n'est alors autre que $(Y_b/\Delta(f_b))$, si $f_b:X_b\to Y_b$ est la restriction de $f$ au-dessus du point g\'en\'erique $b\in B$.
\end{example}

A \'equivalence birationnelle pr\`es, les composantes $\Delta_j$ de $\Delta$ qui sont verticales, c'est-\`a-dire contenues dans une fibre de $g$, peuvent donc \^etre omises. Un tel couple $(X_B/\Delta_B)$ est donc seulement l'un des {\it mod\`eles} de cette clase d'\'equivalence.

\subsection{Courbes orbifoldes}

On va maintenant consid\'erer le cas dans lequel $X_B$ est une courbe sur $k_B$ (et $X$ est donc une surface). Dans ce cas, $\Delta_j$ est une multisection de $g$ (omettant les \'eventuelles composantes verticales). Apr\`es changement de base fini $v:B'\to B$, on peut supposer (et on supposera dans la suite) que $\Delta_j$ est l'image d'une section $p_j:B\to X$ de $g$.

On d\'efinira l'ensemble $M_0\subset B$ des {\it mauvaises places} de $X_B/\Delta_B$ comme \'etant l'ensemble fini des $b\in B$ tels que: ou bien $X_b$ n'est pas lisse, ou bien $X_b$ est lisse, mais la restriction de $g$ \`a $\Delta_B$ n'est pas \'etale au-dessus de $b$.
On choisira aussi arbitrairement un sous-ensemble fini $M\subset B$ qui contient $M_0$.

\begin{example} 
 Soit $g:X\to B$ une courbe sur $k_B$ dont la fibre g\'en\'erique est isomorphe \`a $\Bbb P^1$. A \'equivalence birationnelle pr\`es, $X=\Bbb P^1\times B$, et $g$ est la seconde projection. On notera $\pi:X\to \Bbb P^1$ la premi\`ere projection. 

Soit $\Delta_B=\sum_{j=1}^{j=N} (1-1/m_j).\Delta_j$ un diviseur orbifolde sur $X_B$. On supposera que les $\Delta_j$ sont les images de sections $p_j:B\to X$ \`a $g$.  \`A \'equivalence birationnelle pr\`es, on peut supposer (et on supposera) que les sections $\bar p_j:=\pi\circ p_j:B\to \Bbb P^1$ sont constantes, et prennent les valeurs $0,1$ et $\infty$ respectivement pour $j=1,2$ et $3$.

Dans cette situation, on dira que $X_B/\Delta_B$ est une orbifolde de type $(\Bbb P^1_B/(m_1,m_2,...,m_N))$ si les $m_j$ sont ordonn\'es de telle sorte que $2\leq m_1\leq...\leq m_N$. (On ne suppose pas que les sections $p_j$ soient constantes si $j\geq 4$).

L'ensemble des mauvaises places de $(\Bbb P^1_B/(m_1,...,m_N))$ est donc constitu\'e des points $b$ en lesquels $p_j(b)=p_k(b)$ pour un $k\neq j$. Si $N\leq 3$, cet ensemble est vide.
\end{example}

\begin{definition} Soit $X_B/\Delta_B$ et $X'/\Delta'_B$ des orbifoldes support\'ees par des courbes $X_B$ et $X'_B$ sur $k_B$. Un morphisme d'orbifoldes $g:X_B/\Delta_B\to X'_B/\Delta'_B$ est une application m\'eromorphe surjective $g:X_B\to X'_B$ au-dessus de $k_B$, qui induit un morphisme $g_b:X_b/\Delta_b\to X'_b/\Delta'_b$, pour $b\in B$ g\'en'erique. 
Ce morphisme est \'etale si sa restriction au-dessus de $b$ g\'en\'erique est \'etale.
\end{definition}

\begin{example}  Soit $g:X\to B$ une $k_B$-courbe elliptique (ie: telle que les $X_b$ lisses soient des courbes elliptiques), munie d'une section $s:B\to X$ prise pour \'el\'ement neutre de la loi de groupe dans les fibres lisses de $g$. Il existe alors un morphisme \'etale (au sens orbifolde pr\'ec\'edent) de degr\'e deux $u:X_B\to \Bbb P^1_B/(2,2,2,2)$ qui ramifie (g\'eom\'etriquement) exactement aux points de $2$-torsion de $X/B$. Ici l'orbifolde $X_B$ est l'orbifolde $(X_B/\emptyset)$. Ce m\^eme morphisme induit aussi un morphisme \'etale $u^+:X_B/(1/2.s(B))\to (\Bbb P^1_B/(2,2,2,4))$ lorsque $X_B$ est munie du diviseur orbifolde $\Delta:=(1/2).s(B)$.

\end{example}

\subsection{Points $k_B$-rationnels d'une courbe orbifolde}

Soit $X_B/\Delta_B$ une orbifolde support\'ee par une courbe $X_B$ d\'efinie sur $k_B$, et $M\subset B$ un ensemble fini contenant les mauvaises places de l'orbifolde $X_B/\Delta_B$. On suppose que les composantes $\Delta_j$ de $\Delta_B$ sont les images de sections $p_j:B\to X$ \`a $g:X\to B$.

\begin{definition}\label{kratcf} Soit $(X/\Delta)(k_B,M)\subset X(k_B)$ l'ensemble des $s:B\to X$, sections de $g:X\to B$, telles que, pour tout $b\in B$, si $b\notin M$, et si $s(b)=p_j(b)$, alors $s$ et $p_j$ sont tangentes \`a l'ordre $m_j$ {\bf au moins} en $b$.
\end{definition}

Remarquons que deux mod\`eles birationnels de $X_B/\Delta_B$ ont les m\^emes points $k_B$-rationnels (au sens de la d\'efinition \ref{kratcf} pr\'ec\'edente) si $M$ est assez grand. Et aussi que si $u:X_B/\Delta_B\to X'_B/\Delta'_B$ est un morphisme d'orbifoldes sur $k_B$, alors $u((X_B/\Delta_B)(k_B,M))\subset (X'_B/\Delta'_B)(k_B,M)$, pour tout $M$ assez grand.

Cette d\'efinition provient de la:

\begin{proposition} \label{fctcf}Soit $f:Z_B\to X_B$ une fibration de la $k_B$-vari\'et\'e $Z_B$ sur la $k_B$-courbe $X_B$. Soit $X_B/\Delta(f)$ la base orbifolde de $f$. Alors $f(Z(k_B))\subset (X_B/\Delta(f))(k_B,M)$, pour tout $M\subset B$ assez grand.
\end{proposition}

{\bf D\'emonstration:} Par hypoth\`ese, $f^*(\Delta_j)\geq m_j.f^{-1}(\Delta_j)$. Donc $(f\circ s)^*(\Delta_j)=s^*(f^*(\Delta_j))\geq m_j(f\circ s)^{-1}(\Delta_j)$, ceci pour tout $j\in J$. (On note $A\geq B$ si le diviseur $A-B$ est effectif). D'o\`u l'assertion si $M$ contient toutes les mauvaises places de $\Delta(f)$ et les points de $b$ au-dessus desquels $f$ n'est pas partout d\'efinie.

\begin{remark} Si l'on avait consid\'er\'e $(X_B/\Delta^*(f))$, on aurait obtenu la condition plus forte: $f\circ s$ et $p_j$ sont tangentes en $b$ \`a un ordre {\bf divisible par} $m_j$. C'est le point de d\'epart de [D97].
\end{remark}

\subsection{Mordell orbifolde sur un corps de fonctions}

\begin{theorem} \label{mocf} Soit $(X_B/\Delta_B)$ une orbifolde d\'efinie sur $k_B$. Si $\kappa(X_B\Delta_B)=1$, alors pour tout $M\subset B$, $(X/\Delta)(k_B,M)_{nc}$ est fini. 

On d\'esigne, si $(X_B/\Delta_B)=(F\times B)/(\Delta\times F)$  est triviale, par $(X/\Delta)(k_B,M)_{nc}$ le sous-ensemble de $(X/\Delta)(k_B,M)$ constitu\'e des sections $s$ qui sont non-constantes, c'est-\`a-dire telles que $\bar s:=\pi\circ s:B\to F$ ne soit pas constante, $\pi:X\to F$ \'etant la premi\`ere projection. On peut v\'erifier que la finitude de $(X/\Delta)(k_B,M)_{nc}$ est ind\'ependante des mod\`eles birationnels et des trivialisations choisies.
\end{theorem}

\begin{remark} La d\'emonstration donn\'ee ici n'est pas effective lorsque $X_b$ est rationnelle. Dans ce cas (o\`u $X_b\cong \Bbb P^1$), on peut facilement (tout comme dans le cas classique o\`u $g(X_b)\geq 2$) rendre effective la borne sur la hauteur. Mais la finitude du nombre de points de hauteur donn\'ee n'est pas effective (voir la seconde \'etape de la d\'emonstration).
\end{remark}

{\bf D\'emonstration:} Soit $g(X_b)$ le genre d'une fibre g\'en\'erique de $g:X\to B$. Si $g(X_b)\geq 2$, l'assertion r\'esulte de \ref{comp} et de r\'esultats classiques ([M63], [G65],[P68]). Si $g(X_b)=1$, il suffit, d'apr\`es \ref{comp} de traiter le cas o\`u $\Delta_B=(1/2).B'$, o\`u $B'$ est l'image d'une section $s:B\to X$. Dans ce cas, il existe un rev\^etement double \'etale au sens des orbifoldes $u:(X_B/((1/2).B')\to (\Bbb P^1_B/(2,2,2,4))$. Par la remarque suivant \ref{kratcf}, on en d\'eduit la finitude de $(X_B/\Delta_B)(k_B,M)$ si l'on \'etablit celle de $(\Bbb P^1_B/(2,2,2,3))(k_B,M)$, ce qui sera fait ci-dessous. Nous sommes donc ramen\'es \`a traiter le cas o\`u $X_b\cong\Bbb P^1$.

\

\begin{notation} \label{conv} Lorsque $X_b\cong \Bbb P^1$, il suffit, toujours d'apr\`es la remarque suivant \ref{kratcf} et \ref{comp}, de traiter les cas suivants: $N=3$, $N=4$ et $\Delta_B=(2,2,2,3)$, et enfin: $N=5$ et $\Delta_B=(2,2,2,2,2)$, o\`u $N$ est le nombre de points de $\Delta_b$, pour $b\in B$ g\'en\'erique.

On supposera donc dans la suite que $X_B=\Bbb P^1\times B$, que $g:X\to B$ (resp. $\pi:X\to \Bbb P^1$) est la premi\`ere (resp. la seconde) projection, que les composantes de $\Delta_B$ sont des sections $p_j:B\to X$ pour $j=1,2,...,N$, et que les fonctions $\bar p_j:=\pi\circ p_j:B\to \Bbb P^1$ pour $j=1,2,3$ sont constantes et prennent respectivement les valeurs $0,1$ et $\infty$. Pour $j=1,2,3$, on notera alors $\Delta_0,\Delta_1$ et $\Delta_{\infty}$ respectivement les images des sections $p_j$.

Si $N=3$, on a donc: $\Delta_B=(1-1/m_0).\Delta_0+(1-1/m_1).\Delta_1+(1-1/m_{\infty}).\Delta_{\infty}$, avec: $1/m_0+1/m_1+1/m_{\infty}<1$.

Si $N=4$, on a donc: $\Delta_B=(1/2).(\Delta_1+\Delta_{\infty}+\Delta_p)+(2/3).\Delta_0$, en notant $\Delta_p$ l'image de $p_4$ (la section non n\'ecessairement constante).

Si $N=5$, on a donc $\Delta_B=(1/2).(\Delta_0+\Delta_1+\Delta_{\infty}+\Delta_p+\Delta_q)$, o\`u $\Delta_p$ et $\Delta_q$ d\'esignent les images des sections non n\'ecessairement constantes $p_4$ et $p_5$.
\end{notation}

Le th\'eor\`eme \ref{mocf} est d\'emontr\'e en deux \'etapes: la premi\`ere (lemme \ref{degborn} ci-dessous) est la majoration du degr\'e de $\bar s:=\pi\circ s:B\to \Bbb P^1$, si $s:B\to X\in (X_B/\Delta_B)(k_B,M)$. La seconde (lemme \ref{discr} ci-dessous) est la d\'emonstration du fait que les \'el\'ements de $(X/\Delta)(k_B,M)_{nc}$ forment, vues dans la vari\'et\'e de Chow de $X$, un ensemble discret.

\begin{lemma}\label{degborn} Dans la situation d\'ecrite dans \ref{conv} ci-dessus, il existe un entier $A$ tel que pour tout $s:B\to X\in (X/\Delta)(k_B,M)_{nc}$, le degr\'e de l'application $\bar s:=\pi\circ s:B\to \Bbb P^1$ soit major\'e par $A$.
\end{lemma}

\begin{remark} L'entier $A$ construit ci-dessous d\'epend a priori de $B$, $\Delta_B$ et du cardinal de $M$. Avec un peu plus de travail, on pourrait le faire d\'ependre de mani\`ere effective seulement de $g(B)$, du degr\'e des applications $p_j$, $j=1,\dots,N$, et du cardinal de $M$ seulement.
\end{remark}

{\bf D\'emonstration:} On va distinguer les trois cas: 

a. $N=3$,

b. $N=4$ et $\Delta_B$ de type $(2,2,2,3)$, 

c.  $N=5$ et $\Delta_B$ de type $(2,2,2,2,2)$.

Le premier cas (a) est ext\^emement simple (mais c'est aussi un cas peu int\'er\'essant dans la version corps de fonctions. Ce cas est par contre central dans la version arithm\'etique, gr\^ace au th\'eor\`eme de Bielyi): soit $y:=u/v$ la coordonn\'ee lin\'eaire globale sur $\Bbb P^1$, de coordonn\'ees homog\`enes $[u:v]$. Soit $dy$ la $1$-forme diff\'erentielle m\'eromorphe associ\'ee sur $\Bbb P^1$. On identifiera, pour simplifier les notations, $dy$ et $\pi^*(dy)$, image r\'eciproque sur $X:=\Bbb P^1\times B$ de $dy$ par $\pi:X\to \Bbb P^1$.

Soit alors $s:B\to X\in (X/\Delta)(k_B,M)_{nc}$. Donc $\bar s:B\to \Bbb P^1$ est non-constante, et ramifie \`a l'ordre $m_0$ au moins en chaque point $b\notin M$ tel que $\bar s(b)=0$. On a une hypoth\`ese similaire pour les points $b\notin M$ tels que $\bar s(b)=1$ ou $\infty$. Soit $d>0$ le degr\'e de l'application $\bar s:B\to \Bbb P^1$, c'est-\`a-dire le nombre de points d'une fibre g\'en\'erique. 

Par la formule d'Hurwitz, $2g(B)-2=-2d+\sum _{b\in B}(r(b)-1)$, o\`u $r(b)$ est l'ordre de ramification de $\bar s$ en $b\in B$. Donc $2g(B)-2\geq \sum_{b\in B'} (r(b)-1)$, o\`u $B'$ est l'ensemble (fini) des $b$ tels que $\bar s(b)=0,1,\infty$. Donc $2g(B)-2\geq [\sum_{b\in B'} r(b)]-[\sum_{b\in B'} (1)]$. Comme $\sum_{b\in B'} r(b)=3d$, et que $\sum_{b\in B'} 1\leq [\sum_{b\in M}1]+[\sum_{b\in B'} (r(b)/m_{\bar s(b)})]\leq \vert M\vert + d.(1/m_0+1/m_1+1/m_{\infty})]$, on obtient: $2g(B)-2\geq -2d+3d-(1/m_0+1/m_1+1/m_{\infty}).d-\vert M\vert$, notant $\vert M\vert$ le cardinal de $M$. 

Puisque $1-(1/m_0+1/m_1+1/m_{\infty}):=\epsilon>0$, on en d\'eduit que: $d\leq (2g(B)-2+\vert M\vert)/\epsilon\leq 42.(2g(B)-2+\vert M\vert)$, car $\epsilon\geq 1/42$ pour toute orbifolde de type g\'en\'eral de support $\Bbb P^1$. Ceci \'etablit donc le lemme lorsque $N=3$. Remarquons qu'en fait cet argument s'applique sans changement si $N$ est arbitraire, et si les composantes de $\Delta_B$ sont des graphes d'applications constantes de $B$ dans $\Bbb P^1$.

\

Nous traitons maintenant le second cas (b) ci-dessus\footnote{En consid\'erant implicitement $(\Omega^1_{(X/\Delta)})^{\otimes m}$, pour $m$ ad\'equat. Voir [C01] pour cette notion.}, en construisant une pluri-forme diff\'erentielle m\'eromorphe $w:=(\sum_{k=0}^{k=6}P_k(x,y)dy^{\otimes 6-k}\otimes dx^{\otimes k})/(y^4.(y-1)^3.(y-p(x))^3)$, o\`u $dx$ est l'image r\'eciproque sur $X$ par $g$ d'une $1$-forme diff\'erentielle m\'eromorphe non-nulle sur $B$ (holomorphe si $g(B)\geq 0$), et o\`u les $P_k(x,y)$ sont des fonctions m\'eromorphes sur $X$, \'egales pour $x\in B$ g\'en\'erique \`a des polyn\^omes en $y$.

On va construire $w$ de telle sorte que pour toute section locale holomorphe $s:U\to X$ de $g$, d\'efinie sur un voisinage ouvert analytique $U$ de $b\in B$, les propri\'et\'es suivantes soient satisfaites,  si $b\notin M'$, $M'$ un sous-ensemble fini ad\'equat de $B$ (d\'ependant de $w $ a priori):

{\bf C1.} Si $s(b)=1,\infty$, ou $p(b)$, et si $s(U)$ est tangente \`a $\Delta_1,\Delta_{\infty},\Delta_p$ respectivement, alors $s^*(w)$ est holomorphe en $b$.

{\bf C2.} Si $s(b)=0$, et si $s(U)$ a en $s(b)$ un contact d'ordre $r\geq 3$ avec $\Delta_0$, alors $s^*(w)$ est holomorphe et s'annule en $b$ \`a l'ordre $r/3$ au moins.

L'existence d'une telle forme $w$ est une cons\'equence imm\'ediate du lemme suivant, qui fournit l'existence au voisinage d'une fibre $\Bbb P^1\times \lbrace b\rbrace$ de $g$:

\begin{lemma}\label{w6} Soit $p:\Bbb D\to \Bbb C$ une fonction holomorphe, $\Bbb D$ \'etant  le disque unit\'e dans $\Bbb C$. Il existe des polyn\^omes non nuls $P_k(y),k=0,1,...,6$, tels que si $w:=(\sum_{k=0}^{k=6}P_k(x,y)dy^{\otimes 6-k}\otimes dx^{\otimes k})/(y^4.(y-1)^3.(y-p(x))^3)$, o\`u les $P_k(x,y)$ sont des fonctions m\'eromorphes sur $Z:=\Bbb D\times \Bbb P^1$ d\'ependant holomorphiquement de $x\in \Bbb D$, et qui sont des polyn\^omes en $y$ pour chaque $x$ fix\'e, et tels que: $P_k(0,y)=P_k(y)$, alors la propri\'et\'e suivante est satisfaite: $s^*(w)$ est holomorphe en $x=0$
si $s:\Bbb D\to \Bbb P^1$ est holomorphe, dans les trois cas suivants:

0. $s(x)$ est tangente en $x=0$ \`a la section $\infty$.

1. $s(x)$ est tangente en $x=0$ \`a la section constante de valeur $1$.

2. $s(0)=p$ et $s'(0)=p'$.

3. $s^{(h)}(0)=0$ pour $h=0,1,2$; dans ce cas $s^*(w)$ s'annule en $x=0$. (On a not\'e $s^{(h)}$ la d\'eriv\'ee $h$-i\`eme de $h$). Plus g\'en\'eralement: si $s$ s'annule en $x=0$ \`a l'ordre $r \geq 3$, alors $s^*(w)$ s'annule en $x=0$ \`a l'ordre $(r-1)/2\geq r/3$ au moins.

\end{lemma}

{\bf D\'emonstration:} Il s'agit essentiellement d'un lemme d'alg\`ebre lin\'eaire. La condition 0. est satisfaite si $w(x,u(x)/x^s)$ est holomorphe en $x=0$ pour $s\geq 2$ et $u$ holomorphe avec $u(0)\neq 0$ . Supposons $s=2$. Soit $d_k$ le degr\'e du polyn\^ome $P_k(y)$ (\`a d\'eterminer). La condition est satisfaite pour $s=2$ si, pour chaque $k$, on a: 
$x^{2(4+3+3)}/x^{2d_k+3(6-k)}$ est holomorphe. Cette condition est satisfaite si $d_k\leq 1+3k/2$. On obtient donc la condition 0 si $d_k\leq 1,2,4,5,7,8,10$ pour $k=0,1,2,3,4,5,6$ respectivement. L'espace vectoriel complexe engendr\'e par les coefficients des polyn\^omes $P_k$ est donc de dimension $(1+2+4+5+7+8+10+7)=44$. On v\'erifie facilement que si ces conditions sont satisfaites pour $s=2$, elles le sont aussi pour $s\geq 3$.
 
 La condition $1$ est satisfaite si $w(x,1+x^s.u(x))$ est holomorphe en $0$ pour $s\geq 2$. Ici encore, la consid\'eration du cas $s=2$ suffit. Il suffit que pour tout $k$, $(P_k(x,1+x^2.u(x)).x^{(6-k)})/x^6)$ soit holomorphe en $0$. Ceci est vrai pourvu que les d\'eriv\'ees relatives \`a la variable $y$ satisfassent \`a: $P_k^{(h)}(1)=0$ si $h<k/2$. Les d\'eriv\'ees (relatives \`a la variable $y$) de $P_k$ en $y=1$ doivent donc s'annuler jusqu'\`a l'ordre $a_k$ inclus, avec $a_k=0,0,1,1,2,2$ si $k=1,2,3, 4,5,6$ respectivement. (Il n'y a pas de condition sur $P_0$). 
 
 La condition $2$ est enti\`ement similaire \`a la pr\'ec\'edente, mais appliqu\'ee au point $y=p$, et aux polyn\^omes $Q_k(y)$ d\'efinis comme suit: on doit v\'erifier cette fois-ci que $w(x,p(x)+x^2.u(x))$ est holomorphe en $0$. 
 
 Ecrivons: $w=(\sum_{k=0}^{k=6}Q_k(x,y)dz^{\otimes 6-k}\otimes dx^{\otimes k})/(y^4.(y-1)^3.(y-p(x))^3)$, o\`u $dz=dy-p'(x)dx$. Puisque $dy=dz+p'(x)dx$, on obtient: 
 
 $Q_k=\sum_{m=0}^{m=k}(\matrix{6-k+m\cr m})P_{k-m}.p'(x)^m$. Par l'argument utilis\'e pour la condition $1$, on voit donc que la condition $2$ est satisfaite si $Q_k^{(h)}(p)=0$ pour $h<k/2$, et pour tout $k=1,2,3,4,5,6$.
 
 Le m\^eme argument que pour la condition $1$ montre que la condition $3$ est satisfaite si $P_k^{(h)}(0)=0$ pour $h\leq 2k/3:=b_k$. On a donc: $b_k=0,0,1,2,2,3,4$ si $k=0,1,2,3,4,5,6$ respectivement. On v\'erifie aussi que si $s$ s'annule en $x=0$ \`a l'ordre $r\geq 3$, alors $s^*(w)$ s'annule en $x=0$ \`a l'ordre $(b_k-k+3)(r-1)+(b_k-3)$ au moins. Cette quantit\'e est au moins \'egale \`a $(r-1)/2$ si $r\geq 3$, par v\'erification directe. 
 
 Les polyn\^omes cherch\'es sont les solutions dans la somme directe $W$ des espaces vectoriels de polyn\^omes complexes de degr\'es $d_k$ au plus pour $k=0,...,6$ d'\' equations lin\'eaires homog\`enes. Or $W$ est de dimension $44$. Les conditions $1$ et $2$ n\'ecessitent chacune la satisfaction de $\sum_{k=1}^{k=6}(a_k+1)=12$ conditions (lin\'eaires homog\`enes). La condition $3$ necessite la satisfaction de $\sum_{k=0}^{k=6}(b_k+1)=19$ conditions. L'espace vectoriel des solutions est donc de dimension au moins $44-(12+12+19)=1$. D'o\`u le lemme.
 
 La construction d'une pluri-forme m\'eromorphe globale $w$ satisfaisant (entre autres) les conditions {\bf C1.} et {\bf C2.}  ci-dessus est maintenant imm\'ediate: On consid\`ere sur $X=\Bbb P^1\times B$ le faisceau:
 
  $\cal F$$:=\bigoplus_{k=0}^{k=6}[(dy^{(6-k)}\otimes(g^*(\Omega^1_B)^{\otimes k})/(y^4.(y-1)^3.(y-p)^3)]\otimes [\pi^*(\cal O$$_{\Bbb P^1}(d_k)]$. Le faisceau (sur $B$) $\cal G$$\subset g_*(\cal F)$ egendr\'e par les germes de sections dont les \'el\'ements satisfont en chaque point $b$ de $B$ aux conditions d\'ecrites pr\'ec\'edemmment ($P_k^{(h)}(b,0)=0$ si $h\leq b_k$ pour la condition $3$, par exemple; plus les conditions similaires pour $P_k^{(h)}(b,1)$ et $Q_k^{(h)}(b,p(b))$) est un sous-faisceau coh\'erent de $g_*(\cal F)$ de rang au moins $1$ par le lemme \ref{w6} pr\'ec\'edent. Donc $\cal G$ admet des sections m\'eromorphes globales $\bar w$ non nulles. Une telle section fournit une forme $w$ satisfaisant les conditions voulues. Observons que, par cette construction, $w$ poss\`ede aussi la propri\'et\'e additionnelle {\bf C3.} suivante (en consid\'erant l'ordre maximum $T$ des p\^oles d'une telle section $\bar w$ de $\cal G$):
  
  {\bf C3.} Si $s:B\to X$ est une section de $g$, alors $s^*(w)$ a, en $b$, point arbitraire de $B$, un p\^ole d'ordre au plus $T+4$ (on se ram\`ene en multipliant $w$ par $(x-x(b))^T$, si $x$ est une coordonn\'ee locale sur $B$, au cas o\`u  les coefficients des $P_k$ sont holomorphes en $b$; alors c'est \'evident si $s(b)\neq \infty$; sinon  $s(b)=\infty$, et ceci r\'esulte du calcul fait ci-dessus pour la v\'erification de la condition $0$, en prenant $s=1$). De plus (par le m\^eme argument), si $s$ s'annule en $b$ \`a l'ordre $r\geq 3$, alors la valuation (ordre d'annulation, \'eventuellement n\'egatif) de $s^*(w)$  en $b$ est au moins \'egale \`a $(-T+(r-1)/2)\geq (-T+r/3)$.
 
 Montrons maintenant comment l'existence d'une forme $w$ satisfaisant les conditions ci-dessus entraine le lemme \ref{degborn} dans le cas (b). 
 
 Soit $s:B\to X$ une section de $g$ telle que pour tout $b\notin M$, si $s(b)\in \Delta_j$, alors $s(B)$ est tangente \`a $\Delta_j$ \`a l'ordre $m_j$ au moins, avec $m_0=3$, et $m_j=2$ si $j=1,p,\infty$.
 
 Donc $s^*(w)$ est holomorphe sur $B-M$, et a des p\^oles d'ordre au plus $T+4$ en chaque point de $M^+$, si $M^+$ est la r\'eunion de $M$ et de l'ensemble (fini) des $b\in B$ o\`u l'un des coefficients de l'un des $P_k$ a un p\^ole. De plus, $s^*(w)$ s'annule \`a l'ordre au moins $(r-1)/2\geq r/3$ en chacun des $b\in B$ en lesquels $s(b)$ s'annule \`a l'ordre $r\geq 3$. 
 
 Le degr\'e de $s^*(w)$ est \'egal d'une part \`a $12.(g(B)-1)$, d'autre part \`a $Z-P$, si $Z$ est le nombre de ses z\'eros, et $P$ celui de ses p\^oles, multiplicit\'es comprises. On note $M^+$ la r\'eunion (finie) de $M$ et de l'ensemble des p\^oles de $\bar w$, son cardinal est not\'e $m^+$.
 
 Or $Z\geq \sum_{b\in Q'}(r(b)/3)+\sum_{b\in Q"}(-T+(r(b)/3))$, par la propri\'et\'e {\bf C3.} ci-dessus, si $Q'$ (resp. $Q"$) d\'esigne l'ensemble fini des $b\in B$ tels que $s(b)=0$ et $b\notin M^+$ (resp. tels que $s(b)=0$ et $b\in M^+$). Donc $Z\geq d/3-m^+.T$.
 
 Par ailleurs, $P\leq m^+.(T+4)$, par la condition {\bf C3.} 
 
On en d\'eduit que $12(g(B)-1)\geq d/3-m^+.(T+4)-m^+.T$, c'est-\`a-dire: 

$d\leq 36(g(B)-1)+3m^+.(2T+4)$, et donc le lemme \ref{degborn} dans le cas (b).

\
 
 Nous traitons maintenant le cas (c) de la m\^eme fa\c con. Cette fois-ci, on cherche une forme $w=(\sum_{k=0}^{k=4}P_k(x,y)dy^{\otimes 4-k}\otimes dx^{\otimes k})/(y^2.(y-1)^2.(y-p(x))^2.(y-q(x))^2)$, les diviseurs orbifoldes $\Delta_p$ et $\Delta_q$ \'etant les graphes d'applications $p,q:B\to \Bbb P^1$. 
 
 Les arguments sont tous exactement les m\^emes que dans le cas pr\'ec\'edent (et nous ne les r\'ep\`eterons donc pas), lorsque le lemme suivant est \'etabli:
 
 \begin{lemma}\label{w4} Soit $p,q:\Bbb D\to \Bbb C$ des fonctions holomorphes. 
 
 Il existe des polyn\^omes non nuls $P_k(y),k=0,1,...,4$, tels que si $w:=(\sum_{k=0}^{k=4}P_k(x,y)dy^{\otimes 4-k}\otimes dx^{\otimes k})/(y^2.(y-1)^2.(y-p(x))^2.(y-q(x))^2)$, o\`u les $P_k(x,y)$ sont des fonctions m\'eromorphes sur $Z:=\Bbb D\times \Bbb P^1$ d\'ependant holomorphiquement de $x\in \Bbb D$, et qui sont des polyn\^omes en $y$ pour chaque $x$ fix\'e, et tels que: $P_k(0,y)=P_k(y)$, alors la propri\'et\'e suivante est satisfaite: $s^*(w)$ est holomorphe en $x=0$
si $s:\Bbb D\to \Bbb P^1$ est holomorphe, dans les quatre cas suivants:

0. $s(x)$ est tangente en $x=0$ \`a la section $\infty$.

1. $s(x)$ est tangente en $x=0$ \`a la section constante de valeur $1$.

2. $s(0)=p$ et $s'(0)=p'$.

3. $s(0)=q$ et $s'(0)=q'$

3. $s(0)=s'(0)=0$; dans ce cas $s^*(w)$ s'annule en $x=0$. Plus pr\'ecis\' ement: si $s$ s'annule en $x=0$ \`a l'ordre $r \geq 2$, alors $s^*(w)$ s'annule en $x=0$ \`a l'ordre $r/2$ au moins.
 \end{lemma}
 
 {\bf D\'emonstration:} Elle est identique \`a celle de \ref{w6} dans son principe, en plus simple. 
 
 La condition $0$ se traduit par les conditions similaires: $d_k\leq 2+3k/2$. Donc $d_k\leq 2,3,5,6,8$ pour $k=0,1,2,3,4$. La condition $1$ se traduit par l'annulation des $P_k^{(h)}(1)$ pour $h\leq a_k$, avec $a_k=[(k-1)/2]$. Donc: $a_k=-1,0,0,1,1$ si $k=0,1,2,3,4$. 
 
 Les conditions $2,3$ se traduisent par l'annulation des $Q_k^{(h)}(p)$ et des $R_k^{(h)}(q)$ si $h\leq a_k$, o\`u $Q_k$ et $R_k$ sont les polyn\^omes respectivement d\'efinis par:
 
  $w=(\sum_{k=0}^{k=4}Q_k(x,y)dz^{\otimes 4-k}\otimes dx^{\otimes k})/(y^2.(y-1)^2.(y-p(x))^2.(y-q(x))^2)$ et:
  
  $w=(\sum_{k=0}^{k=4}R_k(x,y)d\zeta^{\otimes 4-k}\otimes dx^{\otimes k})/(y^2.(y-1)^2.(y-p(x))^2.(y-q(x))^2)$, avec:
  
  $dz:=dy-p'(x)dx$ et $d\zeta=dy-q'(x)dx$.
 
 Enfin, la condition $4$ se traduit par l'annulation de $P_k^{(h)}(0)$ pour les $h\leq b_k:=[k/2]$. On a donc: $b_k=0,0,1,1,2$ si $k=0,1,2,3,4$ respectivement.
 
 On voit donc que cette fois-ci, l'espace vectoriel $W$ est de dimension: $2+3+5+6+8+5=29$, tandis que le nombre de conditions (lin\'eaires homog\`enes) est de $3.\sum_k (a_k+1)+\sum_k (b_k+1)=3.(0+1+1+2+2)+(1+1+2+2+3)=18+9=27$. L'espace vectoriel des solutions est donc de dimension au moins $2$.
 
 Ceci ach\`eve la d\'emonstration du lemme \ref{degborn}. 
 
 \
 
 On aborde maintenant la seconde \'etape de la d\'emonstration du th\'eor\`eme \ref{mocf}.
 
 On peut identifier les \'el\'ements $s:B\to X$ de $(X/\Delta)(k_B,M)$ \`a leurs images (r\'eduites), et donc \`a des points de $Chow(X)$, la vari\'et\'e de Chow de $X$. Dans cette identification, $(X/\Delta)(k_B,M)$ est une sous-vari\'et\'e alg\'ebrique ferm\'ee de $Chow(X)$, puisque les conditions qui la d\'efinissent sont alg\'ebriques (les ordres de contact aux points d'intersection non au-dessus de $M$ sont au moins \'egaux aux $m_j$). Le degr\'e des $\bar s$ \'etant born\'e, la finitude sera donc \'etablie si l'on montre que les points de $(X/\Delta)(k_B,M)$, dans cette identification, sont isol\'es.
 
 Il nous suffit, par les arguments usuels du d\'ebut de la d\'emonstration de \ref{mocf}, de d\'emontrer cette propri\'et\'e dans les cas o\`u $X=\Bbb P^1\times B$, et o\`u $\Delta=\sum_{j=1}^{j=N}(1-1/m_j).\Delta_j$, dans les trois cas suivants, o\`u  les $\Delta_j$ sont, pour $j=1,2,3$, les graphes d'applications constantes de $B$ dans $\Bbb P^1$, de valeurs $0,1,\infty$  respectivement:

 a. $N=3$,

b. $N=4$ et $\Delta_B$ de type $(2,2,2,3)$, 

c.  $N=5$ et $\Delta_B$ de type $(2,2,2,2,2)$.

\

\begin{lemma}\label{discr} Dans les trois situations pr\'ec\'edentes, les points de $(X/\Delta)(k_B,M)_{nc}$ sont isol\'es dans $Chow(X)$ si $\kappa(X_B/\Delta_B)=1$ (cette condition est satisfaite dans les cas (b) et (c), elle l'est dans le cas (a) si et seulement si $(1/m_1+1/m_2+1/m_3)<1$).
\end{lemma}

{\bf D\'emonstration:}  Elle est bas\'ee de mani\`ere essentielle sur le lemme \'el\'ementaire suivant:

\begin{lemma}\label{inters} Soit $u=u(t,z):\Bbb D\times \Bbb D\to \Bbb C$, $p=p(z):\Bbb D\to \Bbb C$  des fonctions holomorphes, et $m>0$ un entier. On suppose qu'il existe des fonctions holomorphes $g=g(t):\Bbb D\to \Bbb D$ et $h=h(t,z):\Bbb D\times \Bbb D\to \Bbb C$ telles que: $g(0)=0$, $h(0,0)\neq 0$, et $u(t,z)= p(z)+h(t,z). (z-g(t))^m$. Soit $U$ un voisinage de $0$ dans $\Bbb D$.
Pour tout $t$ assez proche de $0$, l'\'equation: $u(t,z)=u(0,z)$ a, dans $U$ au moins $(m-1)$ solutions distinctes si la fonction $g$ n'est pas constante (\'egale \`a z\'ero, donc), et a au moins $m$ solutions (compt\'ees avec multiplicit\'es) si $g(t)\equiv 0$).
\end{lemma}

\begin{remark}\label{interp} L'interpr\'etation g\'eom\'etrique est la suivante: le graphe $G_t$ de l'application $u_t(z):=u(t,z)$ a un contact d'ordre au moins $m$ avec celui de $p(z)$ au point $x(t):=(g(t), p(g(t)))$. La conclusion est que $G_t$ et $G_0$ ont un nombre d'intersection local (au-dessus de $U$) au moins \'egal \`a $(m-1)$ si le point de contact $x(t)$ est mobile, et au moins $m$ s'il est fixe.
\end{remark}

{\bf D\'emonstration:} On peut supposer que $h(t,z)\neq 0$ partout, quitte \`a restreindre les domaines de d\'efinition. 

Si $g(t)\equiv 0$, l'\'equation \`a r\'esoudre est: $h(t,z).z^m=h(0,z).z^m$. Le r\'esultat est donc \'evident.

Si $g(t)$ n'est pas la fonction nulle, l'\'equation est: $h(t,z).(z-g(t))^m=h(0,z).z^m$. Soit $H(t,z)$ la racine $m$-i\`eme holomorphe de la fonction $h(0,z)/h(t,z)$ qui tend vers $1$ quand $(t,z)$ tend vers $(0,0)$, et $\zeta$ une racine $m$-i\`eme arbitraire de $1$. L'\'equation pr\'ec\'edente est satisfaite si et seulement s'il existe $\zeta$ telle que $z-g(t)=\zeta.H(t,z).z$, c'est-\`a-dire si $z(1-\zeta.H(t,z))=g(t)$. L'\'equation a donc bien $(m-1)$ solutions distinctes qui convergent vers $0$ quand $t$ tend vers $0$.

\

On va maintenant, pour d\'emontrer \ref{discr}, proc\'eder par l'absurde, supposant l'existence dans $(X/\Delta)(k_B,M)$ d'une composante irr\'eductible $T$ de dimension strictement positive, ou encore d'une famille alg\'ebrique $\bar s_t:B\to X$ de sections dont les graphes $S_t$ sont tangents \`a l'ordre $m\geq m_j$ \`a  $\Delta_j$ en chacun des points d'intersection non situ\'es au-dessus de $M$. On choisit un point $0\in T$ g\'en\'erique, de telle sorte que l'ordre des points de contact de $S_0$ avec chacun des $\Delta_j$ au-dessus de chacun des points $b\in B$ (y compris ceux de $M$) soit g\'en\'erique, c'est-\`a-dire minimum, et reste donc le m\^eme pour tous les $t\in T$ voisins de $0\in T$.

Le lemme \ref{inters} peut \^etre appliqu\'e au voisinage de chacun des points $b\in B$ tels que $s(b)=p_j(b)$ pour l'un des $j=1,...,N$, en prenant $u_t=s_t$ et $p:=p_j$. Il montre que le nombre d'intersection local de $S_t$ et de $S_0$ pr\`es de $x(0)=(b,p_j(b))$ est au moins \'egal \`a $(m-1)$ (resp. $m$) si $x(t)\neq x(0)$ (resp. si $x(t)=x(0)$) au voisinage de $0\in T$.

Soit $Q$ l'ensemble des points de $B$ tels que $s(b)=p_j(b)$ pour un $j=1,...,N$ au moins. Soit $Q"$ (resp. $Q'$) l'ensemble des points de $Q$ qui sont (resp. qui ne sont pas) dans $M$ (rappelons que $M$ est l'ensemble des places exclues des conditions de tangence) .

Pour tout $b\in Q$, on note $t_j(b)\geq 0$ l'ordre de contact de $S_0$ avec $\Delta_j$ au-dessus de $b$. On notera qu'il y a un unique $j$ tel que $t_j(b)>0$ si $b\in Q'$, mais qu'il peut y avoir plusieurs tels $j$ si $b\in Q"$ (\`a moins que $N=3$, auquel cas les $\Delta_j$ ne s'intersectent pas non-trivialement). On notera aussi $t(b)$ le plus grand des $t_j(b)$, pour $j=1,...,N$.

On d\'eduit du lemme pr\'ec\'edent \ref{inters}, comptant les nombres d'intersection globaux:

$S_0.S_0=S_t.S_0\geq [ \sum_{b\in Q'}(t(b)-1)]+[\sum_{b\in Q"} t(b)]$.

Par ailleurs, le nombre d'auto-intersection $S_0.S_0=(B'+dF)^2=2d$, si $F$ (resp. $B'$) est la classe de cohomologie d'une fibre de $g:X\to B$ (resp. de l'image d'une section constante de $g$). 

On en d\'eduit que:

(*)  $2d\geq [ \sum_{b\in Q'}(t(b)-1)]+[\sum_{b\in Q"} t(b)]$. 

\

On va d\'eduire une contradiction de cette in\'egalit\'e, dans chacun des trois cas a,b,c ci-dessus.

\

Dans le cas (a), l'argument est, \`a nouveau, tr\`es simple:

On a: $[ \sum_{b\in Q'}(t(b)-1)]+[\sum_{b\in Q"} t(b)]=3.d-\vert Q'\vert$, puisque $\sum_{b\in Q}t(b)=3d$.
Or $\vert Q'\vert \leq (1/m_0+1/m_1+1/m_{\infty}).d$. Donc: l'in\'egalit\'e (*) fournit:

$2d\geq 3d-(1/m_0+1/m_1+1/m_{\infty}).d$. Contradiction, puisque $d>0$ et que $\kappa(X_B/\Delta_B)=1$, par hypoth\`ese, de sorte que $(1/m_0+1/m_1+1/m_{\infty})<1$.

\

Pour d\' emontrer les cas b et c ci-dessus, nous aurons besoin d'un autre lemme.

\begin{lemma}\label{branches} Dans les cas b,c ci-dessus, et avec les notations pr\'ec\'edentes, on a:

$\sum_{b\in Q"}t(b)\geq [\sum_{b\in Q"}\sum_{j=1}^{j=N}(1-1/m_j).t_j(b)]-\sum_{j=1}^{j=N}\delta_j.(1-1/m_j)$, o\`u $\delta_j$ est le degr\'e de l'application $p_j:B\to \Bbb P^1$.
\end{lemma}

Avant de d\'emontrer ce lemme, montrons qu'il entraine \ref{discr}, et ach\`eve donc la d\'emonstration du th\'eor\`eme \ref{mocf}.

On a,  par l'in\'egalit\' e (*): $2d\geq [ \sum_{b\in Q'}(t(b)-1)]+[\sum_{b\in Q"} t(b)]\geq$

$[ \sum_{b\in Q}\sum_{j=1}^{j=N} (1-1/m_j).t_j(b)]+[\sum_{b\in Q'}\sum_{j=1}^{j=N}(1/m_j).t_j(b)]-\vert Q'\vert-\sum_{j=1}^{j=N}\delta_j.(1-1/m_j)$. 

\

Mais $[ \sum_{b\in Q}\sum_{j=1}^{j=N} (1-1/m_j).t_j(b)]=\sum_{j=1}^{j=N} (1-1/m_j).(S_0.\Delta_j)=$

$\sum_{j=1}^{j=N} (1-1/m_j).(d+\delta_j)$, par le calcul cohomologique: 

$(S_0.\Delta_j)=(B'+d.F)(B'+\delta_j.F)=d+\delta_j$.

\

Donc $2d\geq [\sum_{j=1}^{j=N}(1-1/m_j). (d+d_j)]+[\sum_{b\in Q'}(1/m_j).t(b)]-\vert Q'\vert-\sum_{j=1}^{j=N}\delta_j.(1-1/m_j)$.

On a la majoration \'evidente: $\vert Q'\vert\leq [\sum_{b\in Q'}\sum_{j=1}^{j=N}(1/m_j).t_j(b)]$, puisque, si $b\in Q'$, alors $t_j(b)=0$ pour tous les $j$, sauf un seul, pour lequel on a: $t_j(b)\geq m_j$.

\

On en d\'eduit donc que:
$2d\geq [\sum_{j=1}^{j=N}(1-1/m_j)].d$, et une contradiction, puisque $d>0$ et que $[\sum_{j=1}^{j=N}(1-1/m_j)]>2$, ce qui est l'hypoth\`ese $\kappa(X_B/\Delta_B)=1$.

\

{\bf D\'emonstration (du lemme \ref{branches}):}  La d\'emonstration est locale en chaque point $b\in Q"$.  
Supposons que $N=4$. Alors on a, au plus, deux $j$ tels que $t_j(b)>0$, et l'un deux est $j=4$ (correspondant \`a la section non constante $p$, puisque les autres ne s'intersectent pas).

S'il y a un seul (ou aucun) $j$ tel que $t_j(b)>0$, on a: $t(b)\geq \sum_{j=1}^{j=4}t_j(b)$.

Soit $k\neq 4$ le second $j$ tel que $t_j(b)>0$, s'il y a deux des indices $j$ tels que $t_j(b)>0$. Rappelons que l'orbifolde est de type $(2,2,2,3)$, et que les sections $p_j$ d'indices $j=2,3,4$ ont des multiplicit\'es $m_j$ \'egales \`a $2$, tandis que la section $p_1$ (constante de valeur $0$) a une multiplicit\'e $m_1:=3$.

On a alors: $t(b)\geq [\sum_{j=1}^{j=4}(1-1/m_j).t_j(b)]$ si $k\neq 1$, et  $t(b)\geq [\sum_{j=1}^{j=4}(1-1/m_j).t_j(b)-(1/2).\tau(b)]$, o\`u $\tau(b)$ est l'ordre de contact en $(0,b)$ de $\Delta_1$ et de $\Delta_4$ si $k=1$. 

En effet: si $k\neq 1$, alors $(1-1/m_j)=1/2$ pour $j=k,4$, et l'in\'egalit\'e $t(b)\geq \sum_{j=1}^{j=4}(1-1/m_j).t_j(b)$ r\'esulte de ce que $t(b)=max\lbrace t_j(b), j=1,\dots,4\rbrace$, et de ce que $t_j(b)=0$ si $j\neq k,4$.

Si $k=1$, l'in\'egalit\'e $t(b)\geq [\sum_{j=1}^{j=4}(1-1/m_j).t_j(b)-(1/2).\tau(b)]$, r\'esulte de ce que $t_j(b)=0$ si $j\neq 1,4$, et se r\'eduit donc \`a: $t(b)\geq [(1/2)t_4(b)+(2/3)t_1(b)-(1/2).\tau(b)]$. Mais $\tau(b)\geq inf\lbrace t_1(b),t_4(b)\rbrace $ (puisque l'ordre de contact est une valuation). Si $a,b$ sont deux nombres r\'eels positifs de maximum $m$ et de minimum $s$, on a toujours: $m+(1/2).s\geq a/2+2b/3$, et donc la conclusion.

On obtient alors le lemme lorsque $N=4$ en faisant la somme des in\'egalit\'es obtenues sur tous les $b\in Q"$, et en remarquant que les contributions $\sum_{b\in Q"} (1/2).\tau(b)$ sont major\'ees par la moiti\'e du nombre d'intersection de $\Delta_1$ et de $\Delta_4$, qui est \'egal \`a $\delta_4$, et que $\delta_j=0$ si $j\neq 4$.

La d\'emonstration lorsque $N=5$ est similaire.

Les multiplicit\'es $m_j$ attach\'ees aux $5$ composantes $\Delta_j$ sont toutes \'egales \`a $2$.

Si donc $b\in Q"$ est tel que au plus deux des indices $j$ sont tels que $t_j(b)>0$, on a l'in\'egalit\'e:
$t(b)\geq [\sum_{j=1}^{j=5}(1-1/m_j).t_j(b)]$, puisque $(1-1/m_j)=1/2$ et que $t(b)\geq t_j(b)$ pour tout $j$.

Il existe au plus trois des indices $j$ en lesquels $t_j(b)>0$, puisque $4$ des $\Delta_j$ ne s'intersectent pas. Et s'il y a trois tels indices, deux de ces indices sont $j=4$ et $j=5$, puisque les autres $\Delta_j$ ne s'intersectent pas. Soit $k$ le troisi\`eme indice. 

On va montrer que dans ce cas: $t(b)\geq (1/2)(t_k(b)+t_4(b)+t_5(b))-(1/2).\tau(b)$, o\`u $\tau(b)$ est l'ordre de contact en $p_k(b)$ de $\Delta_4$ et de $\Delta_5$. Ceci r\'esulte imm\'ediatement de ce que $t(b)\geq t_j(b),\forall j$, et de ce que $\tau(b)\geq inf\lbrace t_4(b),t_5(b)\rbrace$ (l'ordre de contact d\'efinissant une valuation).

Sommant les in\'egalit\'es obtenues sur les $b\in Q"$, on en d\'eduit que:

$\sum_{b\in Q"}t(b)\geq [\sum_{b\in Q"}\sum_{j=1}^{j=5}(1-1/2).t_j(b)]-[(1/2).\sum_{b\in Q"}\tau(b)]\geq [\sum_{b\in Q"}\sum_{j=1}^{j=5}(1-1/2).t_j(b)]-(1/2).(\delta_4+\delta_5)$, puisque $(\delta_4+\delta_5)$ est le nombre d'intersection $(\Delta_4.\Delta_5)$, \'egal \`a la somme des $\tau(b)$ pour $b\in B$. Ceci ach\`eve la d\'emonstration 
du lemme, et donc celle du th\'eor\`eme \ref{mocf}.
\

\begin{remark}\label{ell} Lorsque $X_B$ est une courbe elliptique sur $k_B$ munie d'un point $k_B$-rationnel (ie: une section $p:B\to X$ de la fibration elliptique $g:X\to B$), le r\'esultat pr\'ec\'edent implique donc la  finitude de l'ensemble des sections $s:B\to X$ dont l'image est tangente \`a $p(B)$ en chacun de leurs points d'intersection (\`a l'exclusion eventuelle des points au-dessus de $M$, fix\'e). Il serait int\'eressant d'avoir un d\' emonstration g\'eom\'etrique directe de ce fait, bas\'ee sur la finitude du rang (sur $\Bbb Z$) du groupe de Mordell-Weil de $X_B$. Une telle d\'emonstration pourrait peut-\^etre s'adapter au cas arithm\'etique. Lorsque le rang du groupe de Mordell-Weil de $X_B$ est $1$, il est tr\`es facile de d\'emontrer le r\'esultat, car $s(B)$ rencontre $p(B)$ en un point $x=s(b)=p(b)$ non au-dessus de $M$, ces deux courbes sont tangentes au point d'intersection $x$. Mais $s(B)$ doit recouper {\it transversalement}, arbitrairement pr\`es de $x$, les multisections (\'etales au-dessus d'un voisinage $U$ de $b$ dans $B$) de $N$-torsion lorsque $N$ est grand. Ce qui montre justement que $Ns(B)$ n'est pas tangente \`a $p(B)$ en leurs points d'intersection proches de $x$. Si $s(B)$ ne rencontre $p(B)$ qu'en des points situ\'es au-dessus de $M$, $N.s(B)$ rencontrera transversalement $p(B)$ en des points non-situ\'es au-dessus de $M$, pour $N$ grand, par l'argument pr\'ec\'edent. Donc l'ensemble des sections consid\'er\'ees est bien fini.
\end{remark} 

\subsection{Finitude pour les surfaces de type g\'en\'eral ayant une fibration de type g\'en\'eral sur une courbe}

\begin{theorem}\label{morcf} Soient $g:X\to B$, $h:C\to B$ et $f:X\to C$ des fibrations telles que $g=h\circ f$, dans lesquelles $X,C,B$ sont projectives complexes lisses et connexes de dimensions respectives $3,2,1$. On suppose que:

1. La fibre g\'en\'erique $X_b$ de $g$ est une surface de type g\'en\'eral.

2. $g$ n'est pas birationnellement isotriviale (ie: $X_b$ n'est pas birationnelle \`a $X_{b'}$ si $b\neq b'$ sont g\'en\'eriques dans $B$).

3. Pour $b\in B$ g\'en\'erique, le restriction $f_b:X_b\to C_b$ de $f$ au-dessus de $b$ est une fibration de type g\'en\'eral.

Alors il existe un sous-ensemble alg\'ebrique strict $D\subset X$ tel que  l'ensemble des sections $s:B\to X$ de $g$ dont l'image n'est pas contenue dans $D$ soit fini.
\end{theorem}

{\bf D\'emonstration:} Dans le langage de \ref{varcf}, on a donc une surface de type g\'en\'eral $X_B$ sur $k_B$, et une fibration de type g\'en\'eral $f_B:X_B\to C_B$ d\'efinie sur $k_B$. L'assertion est l'existence d'un sous-ensemble alg\'ebrique strict $D_B\subset X_B$ tel que l'ensemble $X(k_B)$ des points $k_B$-rationels de $X_B$ qui ne sont pas contenus dans $D_B$ est fini.

Or, d'apr\`es la proposition \ref{fctcf}, $f(X(k_B))\subset (C/\Delta(f))(k_B,M))$, si $M\subset B$ est assez grand fini. D'apr\`es le th\'eor\`eme \ref{mocf}, $(C/\Delta(f))(k_B,M))_{nc}$ est fini. Le th\'eor\`eme \ref{morcf} est donc \'etabli si $(C/\Delta(f))$ n'est pas un produit $(T\times B)/(\Delta\times B)$, o\`u $T$ est une courbe projective connexe, et $\Delta$ un diviseur orbifolde (constant) sur $T$.

On suppose donc d\'esormais que $C=T\times B$, et on note $\pi:C\to T$ la projection sur le premier facteur (la seconde projection \'etant $h$). On notera $X_t\subset X$ la fibre de $\pi\circ f:X\to T$ au-dessus de $t\in T$. La restriction de $g$ \`a $X_t$ d\'efinit donc une fibration $g_t:X_t\to B$ dont les fibres sont celles de $f$ au-dessus de $C_t:=\lbrace t\rbrace\times B$, et sont donc de genre $g\geq 2$, puique par hypoth\`ese $X_b$ est de type g\'en\'eral.

On identifie maintenant $X(k_B)$ au sous-ensemble de $Chow(X)$ param\'etrant les courbes r\'eduites et irr\'eductibles $Z$ de $X$ telles que $X_b.Z=1$. Donc $X(k_B)$ est un ouvert de Zariski de $Chow(X)$, r\'eunion d'un ensemble fini ou d\'enombrable de composantes irr\'eductibles $V_n$. 

Nous noterons $X_n\subset X$ le {\it lieu} de $V_n$, adh\'erence de Zariski de la r\'eunion $X_n^*$ des $s(B)$, pour $s\in V_n$; de sorte que $X_n$ est un sous-ensemble alg\'ebrique ferm\'e irr\'eductible de $X$, dont $X_n^*$ contient un ouvert de Zariski dense.

Pour chacune des $V_n$, deux cas se produisent, simultan\'ement pour tous les $s$ de $V_n$: ou bien $\pi\circ f\circ s:B\to T$ n'est pas constante (ie: $f(s(B)\in (C/\Delta(f))(k_B,M)_{nc}$), ou bien $\pi\circ f\circ s:B\to T$ est constante.

Dans le second cas, $\pi\circ f\circ s(B)$ est l'une des $(C/\Delta(f))(k_B,M)$, qui sont en nombre fini. Tous les $X_n$, pour tous les $n$ poss\'edant cette propri\'et\'e sont donc contenus dans un diviseur $D_{nc}$ de $X$. 

On consid\`ere donc d\'esormais uniquement les $n$ pour lesquels le premier cas se produit, et on note $X(k_B)_c$ l'ensemble des $s\in X(k_B)$ telles que $\pi\circ f\circ s:B\to T$ est constante.

Nous allons maintenant appliquer la version effective de la conjecture de Mordell sur les corps de fonctions \'etablie dans [E-V90] (la premi\`ere version effective est [P68]). Elle affirme que si $u:S\to B$ est une fibration non isotriviale d'une surface projective lisse {\it relativement minimale} $S$ sur une courbe $B$ de genre $q$, et si la fibre g\'en\'erique de $u$ est une courbe de genre $g\geq 2$, alors toute section $s:B\to S$ de $u$ satisfait: $s(B).K_{S/B}\leq 4.(g-1)^2.(q-1+\sigma)$, o\`u $s$ est le nombre de fibres singuli\`eres de $u$, et $K_{S/B}:=K_S-u^*(K_B)$ est le fibr\'e canonique relatif.

Lorsque la fibration $u$ est isotriviale, on a une borne \'el\'ementaire sur $s(B).K_{S/B}\leq 2(q-1)$, pour toute section $s$ de $u$. En effet, $u$ est alors une submersion de fibre $F$. Si $s$ n'est pas isol\'ee dans $S(k_B)$, on a: $s(B).K_{S/B}\leq 0$, par la formule d'adjonction. Si $s$ est isol\'ee, et si $S=F\times B$ est un produit au-dessus de $B$, $K_{S/B}$ est l'image r\'eciproque de $K_F$, la borne ci-dessus est \'evidente. On se ram\`ene \`a ce cas en faisant le changement de base \'etale $v:B'\to B$, o\`u $B':=Aut_B(F\times B,S)$, le groupe relatif des automorphismes de $S$ au-dessus de $B$, dont la fibre au-dessus de $b\in B$ est l'ensemble des isomorphismes de $F$ avec $S_b$. Cette fibre est finie, de cardinal au plus $84(g(F)-1)$, puisque $g(F)\geq 2$.

Soit $T^*\subset T$ l'ensemble (Zariski ouvert non vide) des $t\in T$ pour lesquels $X_t$ est lisse. Remarquons que si $\sigma_t$ est le nombre de fibres singuli\`eres de $g_t:X_t\to B$, la fonction $\sigma_t$ est born\'ee sur $T^*$, puisqu'elle est \'egal au nombre d'intersection de $\lbrace t\rbrace\times B$ avec la courbe $\Sigma \subset C$ constitu\'ee des points au-dessus desquels la fibre de $f$ n'est pas lisse. On notera $\sigma$ ce nombre.

La fin de la d\'emonstration utilise malheureusement des r\'esultats dont la profondeur est sans commune mesure avec le reste du texte.

Nous allons maintenant supposer, quitte \`a faire un changement de base $B'\to B$ fini, que $q:=g(B)\geq 2$. Appliquant la th\'eorie des mod\`eles minimaux au-dessus de $T$ ([K-M92]), nous pouvons supposer (sans changer $X(k_B)$) que $X$ est $\Bbb Q$-factorielle \`a singularit\'es terminales et \`a fibr\'e canonique relatif $K_{X/T}$ nef et {\it vaste}\footnote{traduction de ``big"}. Il existe donc un entier $m>0$ tel que $L:=K_X+mg^*(K_B)+m.(\pi\circ f)^*(H)$ soit nef et vaste sur $X$, si $H$ est ample sur $T$.

Il existe donc [Kaw84] un morphisme birationnel $\psi:X\to Y$ et un fibr\'e ample $L'$ sur $Y$ tels que $L:=K_X+m.(g^*K_B)+(\pi\circ f)^*(H)=\psi^*(L')$. Or, $L.Z\leq 4.(g-1)^2.(q-1+\sigma)+2mq$, pour tout $s\in X(k_B)_c$ (remarquer que $X_t$ reste lisse pour $t\in T$ g\'en\'erique, puisque $X$ est \`a singularit\'es terminales, donc isol\'ees).

On en d\'eduit que la famille $\psi_*(X(k_B))$ est born\'ee. Pour conclure, il suffit donc de montrer qu'aucun des $X_n$ n'est \'egal \`a $X$ entier.

Supposant, par l'absurde que c'est le cas, on obtient une surface projective irr\'eductible $S\subset Chow(X)$ dont l'intersection avec $X(k_B)_c$ est Zariski dense dans $S$, et dont le lieu est $X$. Il existe donc une application rationnelle dominante $\varphi:S\times B\to  X$ d'evaluation des sections,  au-dessus de $B$. Ceci implique, par [Mae83], que $g$ est birationnellement isotriviale, contrairement \`a notre hypoth\`ese.

\

\section{Version Arithm\'etique}

On consid\`ere dans cette section une fibration $f:X\to C$ dans laquelle $X$ (resp. $C$) est une surface (resp. une courbe) projective complexe lisse et connexe. On note encore $F$ une fibre lisse de $f$. On suppose que $X,C, f$ sont d\'efinies sur un corps de nombres $k$.

On a vu que si $X$ et $C$ sont de type g\'en\'eral, alors $X(k')\cap U$ est un ensemble fini, si  $U$ est l'ouvert de Zariski de $X$ r\'eunion des fibres lisses de $f$. On voudrait \'etendre cet \'enonc\'e au cas o\`u $X$ est de type g\'en\'eral, et o\`u $f$ est une fibration de type g\'en\'eral (ie: $\kappa(C/\Delta(f))=1$), avec $g(C)\leq 1$.

Nous allons montrer que cet \'enonc\'e se ram\`ene \`a une version orbifolde de la conjecture de Mordell. Cette r\'eduction est exactement analogue \`a celle utilis\'ee dans le cas des corps de fonctions.

\subsection{Points rationnels orbifolde}

Soit $\cal O$$_k$ l'anneau des entiers de $k$, et $(C/\Delta)$ une orbifolde support\'ee par $C$ et d\'efinie sur $k$, ce qui signifie que $C$ et les points de $C$ formant le support du diviseur $\Delta$ sont d\'efinis sur $k$\footnote{Une d\'efinition plus pr\'ecise, sugg\'er\'ee par P. Eyssidieux, est la suivante: $C$ est definie sur $k$, et pour tout $n$ 
entier le diviseur de Cartier  $\sum_{m_j=n} p_j$est defini sur $k$ (ie 
l'ensemble de $p_j$ pour $m_j=n$ est invariant par le groupe de Galois de 
$k$)}. On suppose choisi un mod\`ele de $(C/\Delta)$ d\'efini sur l'anneau $\cal O$$_{k,M}$ des $M$-entiers de $k$, si $M\subset B$ est un sous-ensemble fini contenant l'ensemble des points de mauvaise r\'eduction de $C/\Delta$, notant $B$ l'ensemble des valuations de $k$. Pour chaque valuation finie $v\in B$, on note $p_v$ l'id\'eal premier de cette valuation.

Soit $\Delta:=\sum_{j=1}^{j=N} (1-1/m_j).p_j$, o\`u les $p_j\in C$ sont dans $C(k)$.

Pour tout $x\in C(k)$, pour toute $v\in B$ finie, et pour tout $j$, on d\'efinit un nombre d'intersection arithm\'etique $(x.p_j)_v$ de $x$ et $p_j$ en $v$ par: $(x.p_j)_v:=max\lbrace m\geq 0\rbrace$, o\`u $m$ d\'ecrit l'ensemble des entiers tels que $x$ et $p_j$ aient m\^eme image dans la r\'eduction de $C$ modulo $p_v^m$. (Voir [D97]).

\begin{definition} \label{rato} Soit $(C/\Delta)(k,M)$ l'ensemble des $x\in C(k)$ tels que pour tout $j=1,...,N$, et tout $v\in B$, finie, $v\notin M$, on ait: $(x.p_j)_v\geq m_j$ si $(x.p_j)_v\neq 0$.
\end{definition}

\begin{remark} \label{comp'} Soit $g:(C/\Delta)\to (C'\Delta')$ un morphisme d'orbifoldes d\'efini sur $k$. Si l'on a des mod\`eles de $g,(C/\Delta)$, et $(C'\Delta')$ d\'efinis sur $\cal O$$_{k,M}$, il est imm\'ediat de v\'erifier que $g((C/\Delta)(k,M))\subset (C'/\Delta')(k,M)$ pour $M$ assez grand.
\end{remark}

L'origine de cette d\'efinition est la:

\begin{proposition} \label{imrat} Soit $f:X\to C$ une fibration, d\'efinie sur $k$, d'une surface $X$ sur une courbe $C$, toutes deux projectives et lisses. Soit $(C/\Delta(f))$ la base orbifolde de $f$, suppos\'ee d\'efinie sur $k$. Si l'on a des mod\`eles de $f, S,C,\Delta(f)$ d\'efinis sur $\cal O$$_{k,M}$, $M$ assez grand, alors $f(X(k))\subset (C/\Delta(f))(k,M)$.
\end{proposition}

{\bf D\'emonstration:} Dans des mod\`eles ad\'equats, $f:X\to C$ est d\'efinie par des polyn\^omes homog\`enes \`a coefficients dans $\cal O$$_{k,M}$. Soit $s\in S(k)$, $v\in B$, $v\notin M$. Supposons que les coordonn\'ees de $f(s)$ et de $p_j$ sont \'egales modulo $p_v$. Par hypoth\`ese $f^*(p_j)=\sum_k m_k.F_k$, o\`u $m_k\geq m_j, \forall k$. On peut donc choisir une carte affine de $C$ et une coordonn\'ee $z$ centr\'ee en $p_j$ telles que $z\circ f=(\Pi_kG_k^{m_k})/G_0$, dans un ouvert affine de $X$ contenant $s$, les $G_k$ \'etant des polyn\^omes \`a coefficients dans $\cal O$$_{k,M}$, et $G_0$ un polyn\^ome tel que $G_0(p_j)$ ne soit pas divisible par $p_v$. Modulo $p_v$, l'un au moins des polyn\^omes $G_k$ s'annule en $s$. D'o\`u l'assertion.

\begin{remark} Dans [D97], qui consid\`ere (implicitement) les multiplicit\'es classiques, la notion de point entier de $(C/\Delta(f))$ est diff\'erente. C'est, en conformit\'e avec le probl\`eme trait\'e, dans la situation de \ref{rato}: $m_j$ divise $(x.p_j)_v$, et non: $(x.p_j)_v\geq m_j$.
\end{remark}

\subsection{Conjecture de Mordell Orbifolde}\label{mo}

C'est la suivante:

\begin{conjecture} \label{mo'} Soit $(C/\Delta)$ l'une des $5$ orbifoldes suivantes, d\'efinie sur un corps de nombres $k$. Alors $(C/\Delta)(k,M)$ est un ensemble fini, pour tout $M$.

1. $\Bbb P^1/(2,3,7)$

2. $\Bbb P^1/(2,4,5)$

3. $\Bbb P^1/(3,3,4)$

4. $\Bbb P^1/(2,2,2,3)$

5. $\Bbb P^1/(2,2,2,2,2)$
\end{conjecture}

\begin{remark} On peut d\'eduire de cette conjecture et de [F83] que pour toute orbifolde de courbe $(C/\Delta)$ d\'efinie sur un corps de nombres $k$, $(C/\Delta)(k,M)$ est fini, pour tout $M$ si $\kappa(C/\Delta)=1$. Ceci r\'esulte imm\'ediatemment de [F83] et de \ref{comp'} si $g(C)\geq 2$. Si $g=0$, c'est une cons\'equence de \ref{mo'}, et de \ref{comp}. Si $g(C)=1$, ceci r\'esulte de \ref{mo'}, de \ref{comp}, et de \ref{etale}. 
\end{remark}

Donnons une cons\'equence imm\'ediate de la conjecture pr\'ec\'edente:

\begin{corollary}\label{lang} Soit $f:S\to C$ une fibration de type g\'en\'eral, d\'efinie sur un corps de nombres $k$, dans laquelle $S$ et $C$ sont projectives complexes lisses et connexes, $S$ une surface de type g\'en\'eral, et $C$ une courbe. Si l'on admet la conjecture \ref{mo'}, alors $S(k)\cap U$ est fini, si $U$ est la r\'eunion des fibres lisses de $S$.

En particulier, il existe, si \ref{mo'} est vraie, des surfaces projectives complexes lisses et simplement connexes d\'efinies sur des corps de nombres et non potentiellement denses.
\end{corollary}

{\bf D\'emonstration:} $f(S(k))\subset (C/\Delta(f))(k,M)$ si $M$ est assez grand, d'apr\`es \ref{imrat}. Comme $\kappa(C/\Delta(f))=1$, par hypoth\`ese, $(C/\Delta(f))(k,M)$ est fini, par la conjecture \ref{mo'}. Par [F83], la restriction de $f$ \`a $S(k)\cap U$ est finie. Donc $S(k)\cap U$ est finie. D'o\`u la premi\`ere assertion.

La seconde assertion est obtenue en appliquant la premi\`ere aux surfaces construites dans \ref{ex}.
(On pourrait d'ailleurs en d\'eduire aussi des exemples potentiellement denses simplement connexes lisses en toutes dimension (en prenant des produits, par exemple, puis des sections hyperplanes)).

\begin{remark} \label{abc}La conjecture \ref{mo'} est, du moins lorsque $N=3$, une cons\'equence imm\'ediate de la conjecture $abc$ (voir [E97], par exemple). Cette observation m'a \'et\'e communiqu\'ee par J.L. Colliot-Th\'el\`ene (et P. Colmez). V\'erifions ceci dans le cas particulier o\`u $k=\Bbb Q$, $N=3$ et $M=\emptyset$.

On peut supposer que $\Delta=(1-1/u).0+(1-1/v).1+(1-1/w).\infty$, o\`u $u,v,w$ sont des entiers au moins \'egaux \`a $2$ et tels que $\epsilon:=1-(1/u+1/v+1/w)>0$. Les points de $(\Bbb P^1/\Delta)(\Bbb Q)$ sont alors les points $x:=a/c\in \Bbb P^1$ tels que $a,c\in \Bbb Z$ sont premiers entre eux et tels que: 
$rad(a)^u$ divise $a$, $rad(c)^w$ divise $c$ et $rad(b)^v$ divise $b$, si $b:=c-a$. On a not\'e $rad(a)$ le produit des nombres premiers qui divisent $a$. 

La condition $rad(a)^u$ divise $a$ signifie donc que chacun de ces nombres premiers apparait dans la d\'ecomposition de $a$ en produit de facteurs premiers avec un exposant au moins \'egal \`a $u$. C'est le cas si, par exemple, $a$ est au signe pr\`es une puissance $u$-i\`eme d'entier, cas trait\'e dans [D97].

La conjecture $abc$ affirme que $rad(abc)=rad(a).rad(b).rad(c)\geq C_t.M^{t}$, pour tout $0<t<1$, avec une constante $C_t>0$, si $M:=max\lbrace \vert a\vert ,\vert b\vert ,\vert c\vert \rbrace$, pour tous $a,b,c$ tels que $a+b=c$ . On a donc, si $x=a/c$ est comme ci-dessus: $M^{(1/u+1/v+1/w)}\geq \vert a\vert ^{1/u}.\vert b\vert ^{1/v}.\vert c\vert ^{1/w}\geq rad(abc)\geq C_t. M^{t}$, pour tout $t<1$. Donc $1\geq C_t.M^{t-(1/u+1/v+1/w)}$. Choisissant $1>t >1/u+1/v+1/w$, on obtient donc \ref{mo'} dans ce cas.
\end{remark}

\

\section{Une fibration de type g\'en\'eral sur une surface simplement connexe}\label{ex}

\

On construit une surface projective lisse et connexe $X$ de type g\'en\'eral {\it simplement connexe} admettant une fibration de type g\'en\'eral sur la droite projective complexe. Ceci montre que, contrairement au cas des multiplicit\'es classiques, il n'y a pas d'obstruction topologique (au niveau du groupe fondamental du moins) \`a l'existence de fibrations de type g\'en\'eral au sens non classique. 

La difficult\'e est de construire des fibrations de base $\bP_1$ ayant une fibre multiple simplement connexe au sens non classique. Bien que l'existence de telles fibres quasiment arbitraires soit \'etablie dans [W74] \footnote{R\'ef\'erence [W74]  communiqu\'ee par E. Amerik} , l'approche suivie dans [W74] (g\'eom\'etrie formelle et th\'eor\`eme d'alg\'ebrisation de Grothendieck) ne permet pas de contr\^oler la base . 

Nous donnons ici une construction explicite mais tr\`es sp\'eciale de telles fibres. Il serait int\'eressant d'avoir des m\'ethodes de construction plus g\'en\'erales que celle propos\'ee ici de telles surfaces simplement connexes ayant une fibration de type g\'en\'eral sur $\bP_1$, ou m\^eme seulement une d\'emonstration constructive des r\'esultats de [W74]\footnote{ Une construction tr\`es simple de fibre multiple par rev\^etement ramifi\'e m'a \'et\'e indiqu\'ee par X. Gang. Cette construction fournit une fibre double avec trois composantes, deux triples et une double. Mais les fibres obtenues ne sont pas simplement connexes.} .

\begin{theorem} Il existe des fibrations de type g\'en\'eral $f:X\to \bP_1$, avec $X$ une surface projective lisse simplement connexe de type g\'en\'eral. On peut choisir $X$ et $f$ d\'efinies sur un corps de nombres $k$. 
\end{theorem}

\begin{remark} Les surfaces construites ci-dessous sont minimales, ont un fibr\'e canonique ample et pour nombres de Chern: $c_1^2=m[(m-1).96-17]$ et $c_2=m[(m-1).48+29]$, o\`u $m\geq 5$ est un entier. Donc $1,66..<(c_1^2/c_2)<2$ (une r\'egion banale de la "g\'eographie" des surfaces de type g\'en\'eral).
\end{remark}

\
La construction des surfaces du th\'eor\`eme pr\'ec\'edent est bas\'ee sur la:

\begin{proposition} {\it Soit $D,L,T\subset \bP_2$ trois droites projectives distinctes et concourantes en un point $a\in \bP_2$. Il existe alors une sextique $S\subset \bP_2$ telle que:

1. $S$ est irr\'eductible et ne passe pas par $a$.

2. $S$ rencontre $D$ en $3$ points distincts en lesquels elle est lisse et tangente \`a $D$.

3. $S$ rencontre $L$ en $3$ points distincts qui sont des points  doubles de $S$.

4. $S$ rencontre $T$ en $2$ points distincts qui sont des points triples de $S$.

Si les droites $D,L,T$ sont d\'efinies sur $\bQ$, les $8$ points pr\'ed\'edents peuvent \^etre choisis sur $\bQ$, et  $S$ peut  \^etre d\'efinie sur un corps de nombres.}
\end{proposition}

\

{\bf D\'emonstration:} On suppose $\bP_2$ rapport\'e aux coordonn\'ees homog\`enes $(U:V:W)$ dans lesquelles les \'equations de $L$ (resp. $T$; resp. $D$) sont: $(U=0)$ (resp. $(V=0)$; resp. $(U-2V=0)$).
On note $C$ la cubique $C:=(D+L+T)$ d'\'equation: $U.V.(U-2V)=0$.

On a donc (par annulation de $H^1(\cal O$$_{\bP_2}(3))$) une suite exacte naturelle d'espaces vectoriels complexes:

$0\to H^0(\cal O$$_{\bP_2}(3))\to H^0(\cal O$$_{\bP_2}(6))\to H^0(\cal O$$_{C}(6))\to 0$.

\

L'espace vectoriel $H^0(\cal O$$_{C}(6))$ s'identifie (par restriction \`a $T,L,D$ respectivement sur l'ouvert affine $W\neq 0$ rapport\'e aux coordonn\'ees $(u:=U/W,v:=V/W)$) aux triplets de polyn\^omes 
$(\bar g(u); \bar h(v); \bar k(v))$ de degr\'es au plus $6$ tels que: $(\bar g(0)=\bar h(0)=\bar k(0))$, et: $(\bar k'(0)=2.\bar g'(0)+\bar k'(0))$.

On peut donc trouver des polyn\^omes $(g,h,k)$ de degr\'es respectifs $(2,3,3)$ tels que les polyn\^omes $(\bar g; \bar h; \bar k):=(g^3;h^2;k^2)$ satisfont les conditions pr\'ec\'edentes, et tels que, de plus: $g(0)=h(0)=k(0)=1$, les z\'eros de chacun des trois polyn\^omes $g^3;h^2;k^2$ \'etant deux-\`a-deux distincts.

Soit alors $S_1(u,v)$ un polyn\^ome de degr\'e $6$  tel que ses restrictions \`a $T,L,D$ soient respectivement $g^3;h^2;k^2$. Ceci exprime que la sextique $S_1$ (d\'efinie par l'\'equation $S_1=0$) coupe $T$ (resp. $L$; resp. $D$) en $2$ (resp. $3$; resp. $3$) points distincts $t_1,t_2$ (resp. $l_1,l_2,l_3$; resp. $d_1,d_2,d_3$) en lesquels elle a un contact d'ordre $3$ (resp. $2$; resp. $2$).

On cherche donc $S$ sous la forme: $S=S_1+E.F$, o\`u $E:=uv(u-2v)$ et $F(u,v)$ est un polyn\^ome (non homog\`ene) de degr\'e au plus $3$.

Pour que les $8=2+3+3$ points d'intersection correspondants aient les multiplicit\'es voulues sur la 
sextique $S$ (d'\'equation $S=0$), il faut et il suffit donc que, notant $F_u, F_{uu}, F_{uv}, etc...$ les d\'eriv\'ees partielles de la fonction $F(u,v)$, on ait:

1. $S_u(t_i)=S_{uv}(t_i)=S_{vv}(t_i)=0$ pour $i=1,2$.

2. $S_u(l_j)=0$ pour $j=1,2,3$.

Compte-tenu du fait que $E(t_i)=E(l_j)=E_u(t_i)=0$ pour tous $i,j$, et de la formule de Leibnitz pour les d\'eriv\'ees partielles d'un produit (ici $E.F)$), ces conditions peuvent \^etre r\'ealis\'ees si les $9$ formes lin\'eaires suivantes:

$F(t_i);F_u(t_i);F_v(t_i); F(l_j)$, pour $i=1,2$ et $j=1,2,3$

 sont lin\'eairement ind\'ependantes sur l'espace vectoriel (de dimension $10$) $H^0(\cal O$$_{\bP_2}(3))$:

Or l'espace vectoriel annulateur de ces formes est constitu\'e des \'equations de cubiques planes ayant un point double en chacun des deux points $t_1,t_2$, et passant par les $3$ points $l_j$. Un argument g\'eom\'etrique imm\'ediat montre que ces \'equations sont de la forme: $\zeta.uv^2$, o\`u $\zeta\in \bC$.

Si $S_2$ est une sextique satisfaisant \`a toutes les conditions pr\'ec\'edentes, et si $\zeta\in \bC$ est g\'en\'erique ad\'equat, la sextique d'\'equation $S:=S_2+\zeta.uv^2$ satisfera aussi ces m\^emes conditions, et sera, de plus, lisse aux points $d_1,d_2,d_3$, car $[uv^2]_u(d_j)\neq 0$, pour $j=1,2,3$.

Pour achever la d\'emonstration de la proposition, il faut encore montrer que le membre g\'en\'erique $S_t$ du syst\`eme lin\' eaire $\Lambda'$ de sextiques engendr\'e par $S$ et $S_0:=D+2L+3T$ est irr\'eductible.

Sinon, $S_t$ est somme de $6$ droites, de $3$ coniques lisses, ou de $2$ cubiques.

Le premier cas est exclu, car une telle droite doit passer par l'un des $t_i$, et par l'un des $l_j$. C'est impossible, puisque $\Lambda'$ n'a pas de composante fixe. Le deuxi\`eme cas est aussi exclu, car une conique composante de $S_t$ devrait \^etre tangente \`a $D$ en l'un des $d_j$, et passer par $2$ des $t_i$ et deux des $l_j$ (avec tangence \`a $T$ ou $L$  si deux de ces points sont confondus). Ceci contredit \`a nouveau le fait que $\Lambda'$ n'a pas de composante fixe.
Le cas o\`u $S_t$ serait r\'eunion de deux cubiques est aussi exclu, car $S_t$ ne serait pas lisse en l'un des points $d_j$ au moins, o\`u les deux cubiques devraient s'intersecter, alors que ce point devrait \^etre un point de tangence de $D$ et de $S_t$.

Les assertions de nature arithm\'etique sont imm\'ediates sur la construction.

\

\begin{corollary}{\it Soit $g:\bP_2\to \bP_2$ le morphisme de degr\'e $4$ d\'efini par: $g(X:Y:Z):=(U:V:W)$, avec $U:=X^2+Y^2;v:=XY;W:=Z^2$. Soit $D$ (resp. $L$; resp. $T$) les droites d'\'equations: $U=2V$ (resp. $U=0$; resp. $V=0$), et $a:=(0:0:1)$. Soit $S$ une sextique d\' equation $S(U:V:W)=0$ satisfaisant les conditions de la proposition pr\'ec\'edente relatives aux droites $D,L,T$. Alors $H=g^*(S)$ est une courbe de degr\'e $12$ telle que:

1. $H$ ne passe pas par $a':=(X=0;Y=0;Z=1)$.

2. $H$ rencontre chacune des trois droites $L',L",D'$ d'\'equations respectives $(Y=i.X); (Y=-i.X);(Y=X)$ en $6$ points distincts qui sont des points doubles de $H$.

3. $H$ rencontre chacune des deux droites $T',T"$ d'\'equations $(Y=0); (X=0)$ en $4$ points distincts qui sont des points triples de $H$.

4. Le membre g\'en\'erique $H_t$ du syst\`eme lin\'eaire $\Lambda$ de courbes de degr\'e $12$ engendr\'e par $H$ et $H':=2.(L'+L"+D')+3(T'+T")$ est irr\'eductible.
}
\end{corollary}

\

{\bf D\'emonstration:} Le rev\^etement $g$ est Galoisien de goupe $G\cong \bZ_2\times \bZ_2$, engendr\'e par $\alpha:(X:Y:Z)\to (-X;-Y;Z)$ et $\sigma: :(X:Y:Z)\to (Y;X;Z)$. De plus, on a:

1. $g(D')=D$ et $g^*(D)=2D'$.

2. $g(L')=g(L")=L$, et $g^*(L)=L'+L"$.

3. $g(T')=g(T")=T$, et $g^*(T)=T'+T"$.

4. $g$ est \'etale sur l'ouvert $Z\neq 0; X^2\neq Y^2$, et ramifie \`a l'ordre $2$ le long de la doite $D'$ priv\'ee de ses deux points $(0:0:1)$ et $(1:1:0)$.

Ces propri\'et\'es permettent de d\'eduire imm\'ediatement le corollaire de la proposition, \`a l'exception de l'irr\'eductibilit\'e de $H_t:=g^*(S_t)$.

\

Si $H_t$ n'est pas irr\'eductible, il est donc r\'eunion soit de 4 cubiques, soit de 2 sextiques irr\'eductibles, qui forment une seule orbite sous l'action de $G$.

Le premier cas (de 4 cubiques) n'est pas possible, car $S_t$ a un point triple en chaque point $t_i$, et $g$ est \'etale au-dessus de ces points. Donc une cubique $\Gamma$ composante de $H_t$ devrait passer par $m=2$ (avec une singularit\'e) ou $m=3$ des $4$ points de $g^{-1}(t_i)$, et aussi par $7$ des autres points de $g^{-1}(R)$, si $R$ est l'ensemble des $8$ points d'intersection de $S_t$ et de $C=L+D+T$, avec les notations de la proposition. Ceci entrainerait que $\Gamma$ est une composante fixe de $\Lambda$, et une contradiction.

Le cas dans lequel $H_t$ serait r\'eunion de $2$ sextiques $S'_t$ et $g.S'_t$, pour $g\in G$ d'ordre $2$ (ind\'ependant de $t$) est \'egalement impossible. Ceci sera \'etabli au cours de la d\'emonstration du corollaire 5.5 ci-dessous.
\

\begin{corollary} {\it Soit $H$ et $H'$ les deux courbes planes de degr\'e $12$ d\'ecrites dans la proposition pr\'ec\'edente, $\Lambda$ le syst\`eme lin\'eaire qu'elles engendrent, et $\psi:\bP_2\to \bP_1$ l'application m\'eromorphe associ\'ee d\'efinie par le quotient (des \'equations) $H/H'$. Soit $b:\bP\to \bP_2$ l'\'eclatement de $\bP_2$ en les $26$ points-base du syst\`eme lin\'eaire $\Lambda$. L'application $\varphi:=\psi\circ b:\bP\to\bP_1$ est holomorphe, \`a fibres connexes. Elle a une unique fibre multiple $\Phi_{\infty}:=\varphi^*(\infty)$, de multiplicit\'e $2$, qui est la transform\'ee stricte de $H'$ par $b$. Cette fibre $\Phi_{\infty}$ est simplement connexe. De plus, $\varphi$ n'a aucune fibre qui soit multiple au sens classique (ie: $m^+(\varphi,y)=1$, pour tout $y\in \bP_1$, avec la notation de la remarque 1.2).
}
\end{corollary}

\

{\bf D\'emonstration:} Seules les assertions: $\varphi$ est holomorphe et n'admet pas de fibre multiple au sens classique ne sont pas cons\'equences imm\'ediates de ce qui pr\'ec\`ede.

On va d\'emontrer la premi\`ere au voisinage de l'un des points d'intersection $t'$ de $T'$ et de $H$, dans des coordonn\'ees locales analytiques $(s,t)$ centr\'ees en ce point, et choisies telles que $T'$ a pour \'equation locale: $t=0$. 

Alors $\psi(s,t):=(a.s^3+b.s^2t+c.st^2+d.t^3+R(s,t))/t^3$, o\`u $R(s,t)$ est une fonction analytique nulle \`a l'ordre $3$ en $(0,0)$, pour des coefficients complexes $a,b,c,d$ ad\'equats.

Le point crucial est: $a\neq 0$. Ceci r\'esulte de ce que le nombre d'intersection de $H$ et de $T'$ est $12$, et est la somme des ordres de contact de $H$ avec $T'$ en chacun des $4$ points $t'_i$. Comme chacun de ces ordres de contact vaut au moins $3$, chacun vaut exactement $3$, et donc $a\neq 0$.

Maintenant le diviseur exceptionnel de l'\'eclatement de $\bP_2$ en $t'$ est recouvert par les deux ouverts de cartes $U_1$ et $U_2$, munis de coordonn\'ees respectives $(\sigma:=s/t,t)$, et $(s,\tau:=t/s)$.

Dans $U_1$, $\psi(s,t)=a\sigma^3+b.\sigma^2+c.\sigma+d+R_1(\sigma,t)$, avec $R_1$ analytique.

Dans $U_2$, $(\psi(s,t))^{-1}=[\tau^3/(a+b.\tau+c.\tau^2+d.\tau^3+s.R_2(s,\tau))]$, avec $R_2$ analytique.

Puisque $a\neq 0$, on en d\'eduit que $\varphi$ est analytique (\`a valeurs dans $\bP_1$) au voisinage de $t'$. 

\

Le cas des points d'intersection de $H$ et $T"$ est identique. Celui des points d'intersection de $H$ avec $L'+L"+D'$ est similaire (en plus simple, car seuls des termes quadratiques interviennent).

Il reste \`a montrer que $\varphi$ n'a pas de fibre multiple au sens classique. Sinon, on aurait pour une valeur ad\'equate de $t\in \bP_1, t\neq\infty $: $H_t=(12/d).\Gamma$, pour $d\neq 12$ un diviseur de $12$, et $\Gamma$ une courbe de degr\'e $d$.

On ne peut pas avoir $d=1,2,3$, car $\Gamma$ doit passser par les $26$ points d'intersection de $H$ et de $H'$. On ne peut pas non plus avoir $d=6$, car la sextique $\Gamma $ devrait passer triplement par chacun des $4$ points d'intersection de $T'$ avec $H$. 

Nous pouvons maintenant montrer que les fibres de $\phi$ sont connexes. Sinon, par l'assertion 4 du corollaire 5.4 pr\'ec\'edent, d\'emontr\'e sauf dans le cas o\`u $H_t$ serait r\'eunion de $2$ sextiques, $\phi$ devrait se Stein-factoriser sous la forme: $\phi=\rho\circ\eta$, avec $\eta:P\to \bP_1$ connexe et $\eta:\bP_1\to \bP_1$ de degr\'e $2$. Alors $\eta$ aurait deux points de ramification double et en ces points, $H_t$ serait une sextique double. Ceci contredit l'argument pr\'ec\'edent, et ach\`eve la d\'emonstration du corollaire.

\

On peut maintenant conclure la d\'emonstration du th\'eor\`eme 4.1:

\begin{corollary}{\it Soit $\varphi:\bP\to Q:=\bP_1$ la fibration pr\'ec\'edente. Soit $g:P:=\bP_1\to Q=\bP_1$ un morphisme fini de degr\'e $m\geq 1$, \'etale au-dessus des points (en nombre fini) de $Q$ au-dessus desquels la fibre de $\varphi$ n'est pas lisse. Soit $X$ la normalis\'ee de $\bP\times _{Q} P$, et $f:X\to P$ le morphisme d\'eduit de $\varphi$ par le changement de base $g$. Soit enfin $u:X\to \bP$ le morphisme fini naturel. Alors:

1. $X$ est une surface lisse, et $f$ est une fibration de type g\'en\'eral si $m\geq 5$.

2. $X$ est simplement connexe (et de type g\'en\'eral si $m\geq 5$).

3. Si la sextique $S$ est d\'efinie sur un corps de nombres, ainsi que $g$, alors $X$ et $f$ sont aussi d\'efinies sur un corps de nombres.
}
\end{corollary}

\

\begin{remark} Si $F$ est une fibre lisse de $\varphi$, on a: $g(F)=13$. En effet, si $E$ (resp. $G$) est le diviseur exceptionnel de $b$ au-dessus des 18 points doubles (resp. 8 points triples) de $H$ consid\'er\'es, alors $K_P=b^*(K_{\bP_2})+E+G$. De plus, si $F$ est une fibre lisse de $\varphi$, on a: 
$F.E=2. 18$, et $F.G=3. 8$ (car chaque composante de $E$ (resp. $G$) se projette doublement (resp. triplement) sur $\bP_1$ par $\varphi$). Donc $K_P.F=2(g(F)-1)=(b^*(K_{\bP_2})+E+G).F=-3.12+2.18+3.8=24$.
\end{remark}

\

{\bf D\'emonstration de 5.6:} Toutes les assertions sont imm\'ediates, sauf peut-\^etre le fait que $f$ soit une fibration de type g\'en\'eral, et que $X$ est simplement connexe. 

Le fait que $f$ soit de type g\'en\'eral si $m\geq 5$ r\'esulte de ce que $f$ a maintenant $m$ fibres multiples (au sens de 1.1) de multiplicit\'e $2$ (qui est la multiplicit\'e de $\Phi_{\infty}$). On conclut \`a l'aide de la remarque 1.4.

La simple connexit\'e de $X$ r\'esulte du lemme standard suivant:

\

\begin{lemma} {Soit $f:X\to C$ une fibration holomorphe propre et surjective sur une courbe. Supposons que $f$ n'ait pas de fibre multiple au sens classique (voir 1.2), et que $f$ ait une fibre simplement connexe $F_0$. Alors: $f_*:\pi_1(X)\to\pi_1(C)$ est un isomorphisme de groupes.}
\end{lemma}

\

{\bf D\'emonstration:} Puisque $f$ n'a pas de fibre multiple au sens classique, le noyau de $f_*$ est l'image de $\pi_1(F)$ dans $\pi_1(X)$ par $j_*$, d\'eduite de l'injection naturelle de $F$ dans $X$, si $F$ est une fibre lisse de $f$. Mais $j$ se factorise par l'injection de $U$ dans $X$, si $U$ est un voisinage ouvert de $F_0$ qui se r\'etracte sur $F_0$, et qui contient $F$, choisie assez proche de $F_0$. Puisque $F_0$ est simplement connexe, il en est de m\^eme de $U$. D'o\`u la conclusion.

\

\section{Bibiographie}

[B-P-V84] W.Barth-C.Peters-A.Van de Ven. Compact Complex Surfaces. Springer Verlag.

[B72] T. Barth.  The  Kobayashi distance induces the standard topology. Proc. Amer. Math. Soc. 35 (1972), 439-441.

[C01] F. Campana. Special varieties, Orbifolds and Classification theory.  math. AG/0110051. et  Ann. Inst. Fourier 54 (2004), 499-665. Pour un expos\'e introductif: math. AG/0402243.

[Ca 04] L. Caporaso. Moduli theory and Arithmetic of Algebraic Varieties. math. AG/0311465.

[C-S-S 97] J.L. Colliot-Th\'el\`ene- A. Skorobogatov- P. Swinnerton-Dwyer. Double fibres and double covers: paucity in rational points. Acta Arithm. 79 (1997), 113-135.

[D97] H. Darmon. Faltings plus epsilon, Wiles plus epsilon, and the generalised Fermat equation. C.R. Math. de l'Acad. des Sciences de la r\'epublique du Canada. 19 (1997), 3-14.

[E97] N. Elkies. ABC implies Mordell. Int. Math. Res. Notes. 7 (1991),99-109.

[E-V90] H. Esnault-E. Viehweg. Effective bounds for semipositive sheaves and the height of points on curves over complex function fields. Comp. Math. 76 (1990), 69-85.

[F83] G. Faltings. Endlichkeitss\" atze f\" ur Abelsche Variet\" aten \" uber Zahlk\" orpern. Inv. Math. 73 (1983), 349-366.

[G65] H. Grauert. Mordells Vermutung \" uber rationale Punkte auf algebraischen Kurven und Funktionenk\" orper. Publ. Math. IHES 25 (1965), 131-149.

[Kaw84] Y. Kawamata. The cone of curves of algebraic varieties. Ann. Math. 119 (1984), 95-110.

[K98] S. Kobayashi. Hyperbolic complex spaces. Springer Verlag (1998).

[K-M92] J.Kolla\`r-S. Mori. Classification of three-dimensional flips. J. Amer. Math. Soc. 5 (1992),533-703.

[Mae83] K. Maehara. A finiteness property of Varieties of general type. Math. Ann. 262, 101-123.

[M63] Y. Manin. rational points of algebraic curves over function fields. Iz. akad. Nauk. SSSR (1963), 1395-1440.

[N87] M. Namba. Branched coverings and algebraic functions. Pitman research notes in Math. Sc. 161. Longman Sc. and Technical. (1987).

[P68] A.N. Parshin. Algebraic curves over function fields I. Math. USSR Izvestija (1968), 1145-1170.

[Siu80] Y.T. Siu. The complex-analyticity of harmonic maps and the strong rigidity of compact k\" ahler manifolds. Ann. Math. 112 (1980), 73-111.

[Vi83]. E. Viehweg. Additivity theorem

[W74] G.Winters. On the existence of certain families of curves. Amer. J. Math. 96 (1974), 215-228.

[Wr77] M. Wright. The Kobayashi pseudometric on manifolds of general type. Trans. AMS. 232 (1977), 357-370.

\end{document}